\documentclass[review,onefignum,onetabnum]{siamart171218}

\usepackage{lipsum}
\usepackage{amsfonts}
\usepackage{graphicx}
\usepackage{epstopdf}
\usepackage{algorithmic}

\usepackage{amssymb}
\usepackage{amsmath}
\usepackage{float}
\usepackage{color}
\graphicspath{{./Images/}}

\ifpdf
  \DeclareGraphicsExtensions{.eps,.pdf,.png,.jpg}
\else
  \DeclareGraphicsExtensions{.eps}
\fi

\newsiamremark{remark}{Remark}
\newsiamremark{hypothesis}{Hypothesis}
\crefname{hypothesis}{Hypothesis}{Hypotheses}
\newsiamthm{claim}{Claim}
\newsiamthm{exmp}{Example}

\headers{FSI for the Classroom: Interpolation!}{N. A. Battista}

\title{Fluid-Structure Interaction for the Classroom: Interpolation, Hearts, and Swimming! \thanks{Submitted to the editors 22/08/2018.
\funding{This work was funded by the NSF OAC-1828163 the Support of Scholarly Activities Grant from TCNJ (The College of New Jersey)}}}

\author{Nicholas A. Battista\thanks{Department of Mathematics and Statistics, The College of New Jersey, 2000 Pennington Road, Ewing Township, NJ 08628 
  (\email{battistn@tcnj.edu}, \url{http://http://battistn.pages.tcnj.edu/}).}}

\usepackage{amsopn}

\ifpdf
\hypersetup{
  pdftitle={Fluid-Structure Interaction for the Classroom: Interpolation, Hearts, and Swimming!},
  pdfauthor={N. A. Battista}
}
\fi

\begin{document}

\maketitle

%
%

\begin{abstract}

While students may find spline interpolation quite digestible, based on their familiarity with continuity of a function and its derivatives, some of its inherent value may be missed when students only see it applied to standard data interpolation exercises. In this paper, we offer alternatives where students can qualitatively and quantitatively witness the resulting dynamical differences when objects are driven through a fluid using different spline interpolation methods. They say, \textit{seeing is believing}; here we showcase the differences between linear and cubic spline interpolation using examples from fluid pumping and aquatic locomotion. Moreover, students can define their own interpolation functions and visualize the dynamics that unfold. To solve the fluid-structure interaction system, the open-source fluid dynamics software \textit{IB2d} is used. In that vein, all simulation codes, analysis scripts, and movies are provided for streamlined use.

\end{abstract}




%
%

%
%
%
%

\section{Introduction}
\label{intro}

Traditionally it is in numerical analysis and scientific computing courses where students are first introduced to the topic of interpolation. It is frequently motivated by posing the seemingly innocent question of, ``If handed $N$ unique data points, $\{x_j,y_j\}_{j=0}^{N}$, can you find a polynomial, $p(x)$, with the property that $p(x_j)=y_j, \forall j=0,1,2,\ldots,N$?" It is customary to accompany this question with a uniqueness theorem that gives a somewhat surprising result for students - that if such a polynomial exists, it must be unique. The proof is even rather elegant \cite{Kincaid:2002,Burden:2014}!

What happens next? Well, surely a discussion of how to construct such a polynomial and alas the standard ways to find such an interpolation polynomial (monomial, Newton, and Lagrange) are derived. This effort, in essence, enforces that students once again see that existence and uniqueness go together, like peas and carrots. 

This may leave the students usually wondering, ``Well, how close is this polynomial to the actual function from which the data was originally sampled?" Not be disappointed, the class dives into estimating the error of such a polynomial, and after seeing a few exploitative examples using uniformly spaced nodes \cite{Runge:1901,Kincaid:2002,Heath:2002,Burden:2014}, and going down the rabbit hole of Chebyshev nodes, students see the corresponding interpolation error and how it can be minimized. 

If that is the best such a polynomial can do in terms of minimizing the error, instructors may encourage their class to contemplate whether there could be any other methods to interpolate the original data given. That is, motivating the students to move beyond constructing a single global polynomial that interpolates the data, but instead interpolating the data point-by-point. This, of course, leads to the introduction of spline interpolation, cubic splines, and/or Bezier curves! Splendid!

Unfortunately, a genuine difficultly for students during this onslaught of interpolation techniques, error analysis, and implementation, is sometimes seeing the practical applications of interpolation. Some possible (surprising) applications for students that may be mentioned include how letters are shaped in typography \cite{Adobe:1997,SVG:2011}, vector graphics and imaging \cite{Vince:2007}, or data and digital signal processing \cite{Unser:1999,Marcus:2016}. However, students generally interested in computational science and modeling may not be captivated or satisfied with these applications. 

We would like to introduce an application of interpolation that unfortunately falls through the cracks for students - the use of interpolation in mathematical modeling, and in particular biological fluid dynamics. Simply stated interpolation can be used to prescribe the motion of objects. The enticing portion - these objects can be immersed within a fluid, where the fluid reacts and moves due to the prescribed motion of said object. 

Not sold, yet? Numerous recent scientific studies have used this exact type of interpolation to successfully prescribe motion, ranging from diverse fields such as heart development \cite{Baird:2014,Lee:2013,Battista:2017}, aquatic locomotion \cite{Hershlag:2011,Alben:2013,Becker:2015}, animal flight \cite{Miller:2009,Ruck:2010,SJones:2015}, organismal feeding and filtering \cite{Hamlet:2012,Nielson:2017,Samson:2017}, and beyond. 

We offer a software alternative that will allow students to test out varying kinds of spline interpolation to prescribe the motion between one or more feature states, within a framework that provides direct practical scientific applications. 

In the remainder of this paper, we will provide three differing examples of how spline interpolation can be used to drive the motion of a structure immersed within a fluid, while also comparing different kinds of spline interpolation, e.g., linear and higher order polynomial (cubic). This will provide students intuition about splines that is not traditionally emphasized in the classroom that can help facilitate greater learning and further curiosity in computational science. 

In Section 2 we motivate the ideas of spline interpolation through the presentation of a moving circular object immersed in a fluid. In Section 3 we introduce how to prescribe motion using a cartoon heart pumping example and provide a stencil for how to create your own example. In Section 4 we move beyond prescribing the motion of individual points to instead interpolate between different material property states of an immersed body, e.g., modeling a structure that has time-dependent curvature, which gives rise to forward locomotion (swimming)! For details regarding the fluid-structure interaction software, see Appendix \ref{IB:Appendix}, or \cite{Battista:2015,BattistaIB2d:2016,BattistaIB2d:2017} for a more detailed overview. All simulations presented here are available on \url{https://github.com/nickabattista/ib2d} and can found in the sub-directory \texttt{IB2d/matIB2d/Examples/Examples$\_$Education/} as well as the Supplementary Materials.

%
%

%
%

\section{Spline Interpolation: Linear vs. Higher Order Polynomials}
\label{Interpolation}
$ $\\

When first introducing splines in numerical analysis, it may fruitful to tell students they have already seen an example of a linear spline in Multivariate Calculus, when parameterizing curves for line integrals. Have them consider two points, \textbf{a} and $\textbf{b}$, $(x_a,y_a)$ and $(x_b,y_b)$, respectively. Students can then parameterize a straight line between the two points in a familiar way:
\begin{equation}
    \label{param_line} (x(t),y(t)) = \textbf{h}_0(t) =\textbf{a} + \frac{t}{t_1}( \textbf{b}-\textbf{a}),
\end{equation}
for $t\in[0,t_1]$. We can see that $\textbf{h}_0(0) = \textbf{a}$ and $\textbf{h}_0(t_1)=\textbf{b}$. Of course, in calculus this is not introduced as a spline and the word interpolation probably doesn't echo off the classroom walls, but that is exactly what this process was - setting up a linear spline interpolant between two points. If we had a third point $\textbf{c}=(x_c,y_c)$, we could construct another linear interpolant between the \textbf{b} and \textbf{c}, 

\begin{equation}
    \label{param_line2} (x(t), y(t)) = \textbf{h}_1(t) = \textbf{b} + \frac{t-t_1}{t_2-t_1}( \textbf{c}-\textbf{b}),
\end{equation}
for $t\in[t_1,t_2]$. We note that $\textbf{h}_1(t_1) = \textbf{b}$ and $\textbf{h}_1(t_2)=\textbf{c}$. The piecewise linear interpolant between all three points could then be written as

\begin{equation}
    \label{param_line_pw} (x(t),y(t)) = \left( \begin{array}{c} \textbf{h}_0(t) \\ \textbf{h}_1(t)   \end{array} \right) = \left\{\begin{array}{ll} 
    \textbf{a} + \frac{t}{t_1}( \textbf{b}-\textbf{a}) & 0\leq t\leq t_1 \ \ \ \ \\ 
    \textbf{b} + \frac{t-t_1}{t_2}( \textbf{c}-\textbf{b}) \ \ \ \ \ & t_1\leq t\leq t_2 \end{array} \right..
\end{equation}

What we have done, although perhaps not emphasized too much in Calculus, is created a method to prescribe the motion of a point, $\textbf{x}$ around the plane in $\mathbb{R}^2$, $$\textbf{a}\rightarrow\textbf{b}\rightarrow\textbf{c}.$$

There is no reason this cannot extend to a larger collection of points! Instead of points \textbf{a}, \textbf{b}, and \textbf{c}, consider the following matrices, where each column contains $N$-$(x,y)$ points, respectively,

\begin{equation}
    \label{matrices} \textbf{A} = \left[ \begin{array}{cc} x_0^a & y_0^a \\ x_1^a & y_1^a \\ \vdots & \vdots \\ x_N^a & y_N^a \end{array} \right] , \  \textbf{B} = \left[ \begin{array}{cc} x_0^b & y_0^b \\ x_1^b & y_1^b \\ \vdots & \vdots \\ x_N^b & y_N^b \end{array} \right] , \ \mbox{ and }\ \textbf{C} = \left[ \begin{array}{cc} x_0^c & y_0^c \\ x_1^c & y_1^c \\ \vdots & \vdots \\ x_N^c & y_N^c \end{array} \right].
\end{equation}

We can write an analogous spline interpolant to (\ref{param_line_pw}) as follows,
\begin{equation}
    \label{param_line_pw2} (\textbf{x}(t), \textbf{y}(t)) = \left( \begin{array}{c} \textbf{H}_0(t) \\ \textbf{H}_1(t)   \end{array} \right) = \left\{\begin{array}{ll} 
    \textbf{A} + \frac{t}{t_1}( \textbf{B}-\textbf{A}) & 0\leq t\leq t_1 \ \ \ \ \\ 
    \textbf{B} + \frac{t-t_1}{t_2}( \textbf{C}-\textbf{B}) \ \ \ \ \ & t_1\leq t\leq t_2 \end{array} \right..
\end{equation}

%
%

\begin{exmp}
\label{ex:circles}

    Consider the circles given by the following $N$ points $\{ x_j^a,y_j^a \}_{j=0}^N, \{ x_j^b,y_j^b \}_{j=0}^N$ and $\{ x_j^c,y_j^c \}_{j=0}^N$. These are illustrated in Figure \ref{Example:3Circles}. 
    \begin{figure}[H]
        \centering
        \includegraphics[width=0.45\textwidth]{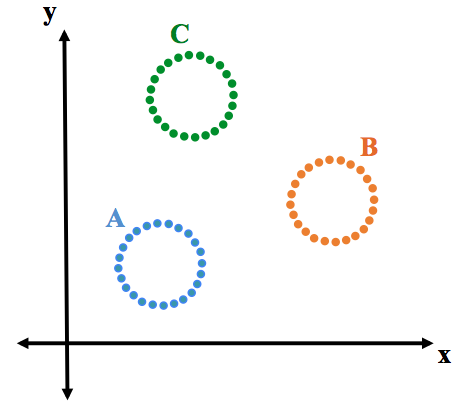}
        \caption{$3$ circles in the $xy$ plane, each composed of $N$ points.}
        \label{Example:3Circles}
    \end{figure}
    
    Next using (\ref{param_line_pw2}), let's prescribe the motion of these circles starting from State \textbf{A} to State \textbf{B} and finally State \textbf{B} to State \textbf{C} for $0\leq t \leq t_2$, with $t_1\in(0,t_2)$ The positions, $(\textbf{x}(t), \textbf{y}(t))$ of these interpolated states are illustrated in Figure \ref{Example:3CirclesInterp}, given by the circle in red. 
    
    \begin{figure}[H]
        \centering
        \includegraphics[width=0.99\textwidth]{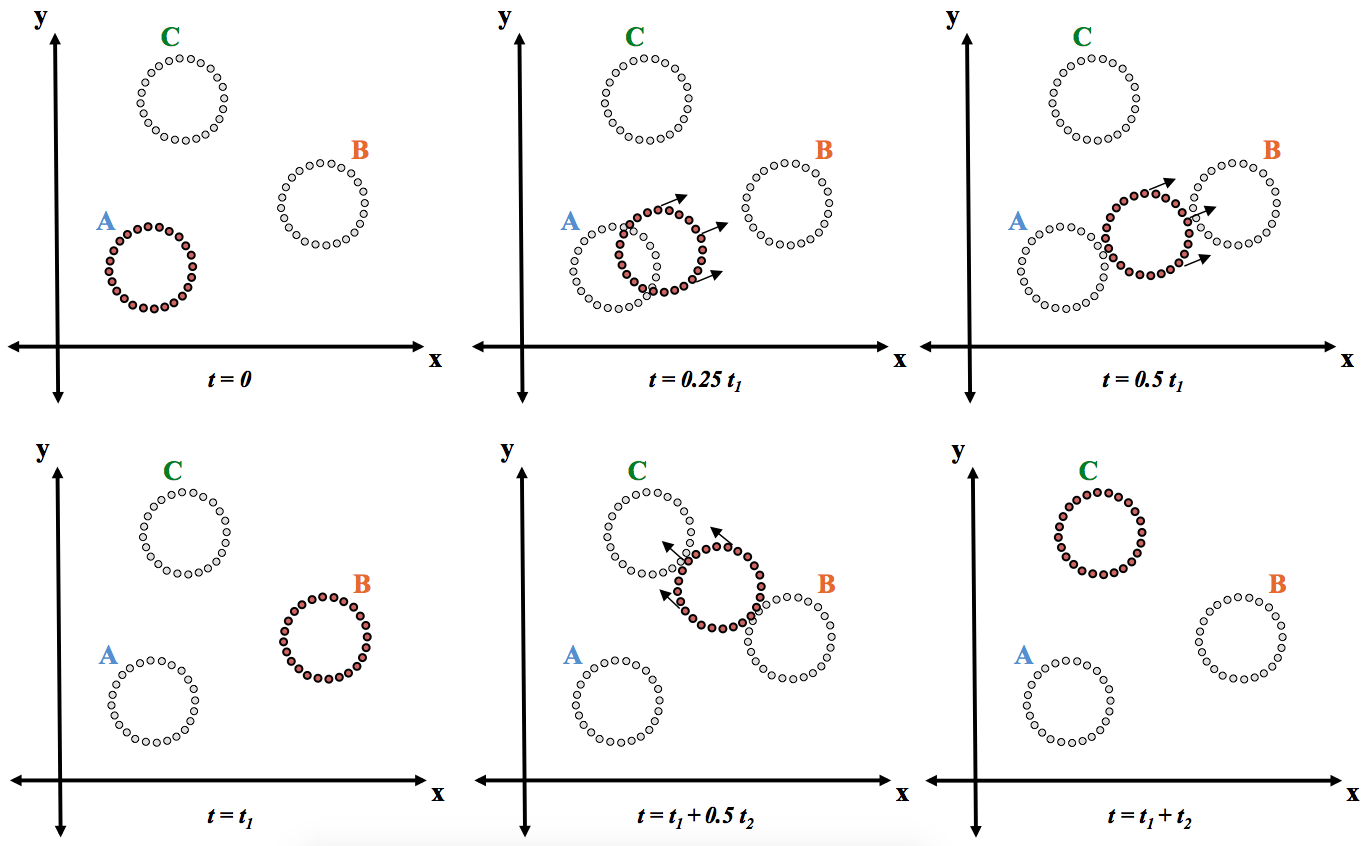}
        \caption{Timeslices of the $N$-point circle moving from State \textbf{A} to State \textbf{B} to State \textbf{C} using the piecewise linear interpolate to prescribe the motion.}
        \label{Example:3CirclesInterp}
    \end{figure}
    
    As mentioned earlier, we could imagine that beyond these circles simply moving around the $xy$-plane in a prescribed fashion, one could envision these objects immersed within a fluid. This is exactly an example found in \textit{IB2d}, e.g.,\\ \texttt{Examples$\_$Education/Interpolation/Moving$\_$Circle/Linear$\_$Interp}. Immersing a circle within a fluid environment and then prescribing its motion will cause the fluid to react, and in turn, move in response. This is shown in Figure \ref{Example:IB_Circles}, where the colormap illustrates the magnitude of the fluid velocity and vector field represents the fluid velocity.
    
    \begin{figure}[H]
        \centering
        \includegraphics[width=0.99\textwidth]{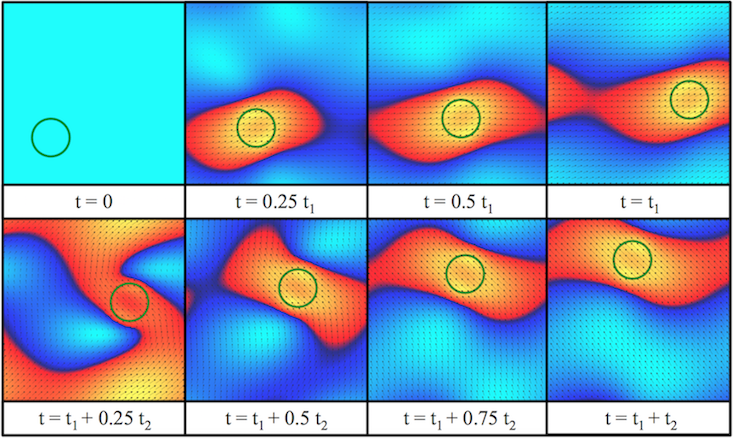}
        \caption{A circle undergoing prescribed motion in a fluid domain, causing the fluid to move in response. The colormap illustrates the magnitude of velocity, while the vector field depicts the fluid velocity itself.}
        \label{Example:IB_Circles}
    \end{figure}
    
    It is evident from Figure \ref{Example:IB_Circles} that the fluid is moving the fastest right nearest to the circle, the immersed object. Students can change the fluid viscosity, $\mu$, or interpolation time-points, $t_1$ or $t_2$, to see how the fluid motion changes. Furthermore, students can plot the simulation as it runs directly within MATLAB, or they can view the data using open-source visualization software, such as VisIt \cite{HPV:VisIt}, which was used to construct Figure \ref{Example:IB_Circles}. Note that this simulation was designed to use a rather unresolved grid, e.g., $32\times32$, for speed so students can watch the movement of the circle unfold directly in MATLAB .
    
    It should be emphasized that while this example only prescribed the motion of a circle, immersed within a fluid, to move between a few predetermined states, this is exactly the kind intepolation that is used in a lot of research applications, as mentioned in Section \ref{intro}. One could imagine constructing a much more complex geometrical entity, such as a heart, fish, or other immersible structure, and prescribe it to move in rather complicated ways in order to test a hypothesis or engineering question!
    
\end{exmp}


From the way the linear interpolant in (\ref{param_line_pw}) and (\ref{param_line_pw2}) was constructed, it should not be a surprise that the interpolant is continuous at all of the interpolation nodes, $\{x_j\}$, that is 
\begin{equation}
    \label{linear_C1} \textbf{h}_{k}(t_{k+1}) = \textbf{h}_{k+1}(t_{k+1}).
\end{equation}

At this stage, students are usually encouraged to consider what happens to the derivatives at the interpolation nodes. Simply differentiating either (\ref{param_line_pw}) or (\ref{param_line_pw2}), one can show that that this linear interpolating scheme does not guarantee continuous derivatives at the nodes. Is this an issue? 

Let's consider the movement of the circle from Example \ref{ex:circles}. When the circle is moving between State A to State B, what happens when $t\approx0$ or $t\approx t_1$? We want to explore how fast the circle moving, its acceleration, and what implications these may have on the circle moving around. There are a couple things to consider: 

\begin{enumerate}
    \item First, we see that going from $t=0$ to $t=\epsilon$, where $\epsilon>0$, that the structure immediately begins to move at a constant speed, the constant speed it will move with between $0\leq t\leq t_1$. This illustrates there is an instantaneous acceleration from not moving to moving at its constant speed.
    \item Second, a similar phenomenon happens as $t\rightarrow t_1$; that is, an instantaneous deceleration from moving at its constant to speed to 0.
    \item Third, if we are testing a hypothesis about the natural world or modeling an engineering device, no such situation occurs where we see such instantaneous accelerations (or decelerations for that matter).
\end{enumerate}

We can encourage students to ask how can we ensure such accelerations do not happen? This can lead to a great discussion on not having enough \textit{degrees of freedom} to enforce continuous derivatives, if we only have piecewise linear interpolating functions. Students may be obliged to try a polynomial of higher degree to interpolate between the positions, such as a quadratic or a cubic. 

Before diving right in, note that the situation we were previously considering had the general linear interpolant
\begin{equation}
    \label{linear_spline} \textbf{h}(t) = \left( \begin{array}{c} \textbf{h}_0(t) \\ \textbf{h}_1(t)   \end{array} \right) = \left\{ \begin{array}{ll} \textbf{a} + (d_0 + d_1t)( \textbf{b}-\textbf{a})   \ \ & 0\leq t\leq t_1 \ \ \ \ \\ \textbf{b} + (d_2 + d_3t)( \textbf{c}-\textbf{b}) \ \ & t_1\leq t\leq t_2 \end{array} \right.,
\end{equation} 

with unknowns, $\{D_j\}_{j=0}^3$. Whether we knew it or not, we constructed (\ref{param_line_pw}) and (\ref{param_line_pw2}) using the following continuity conditions to find the unknown coefficients:


\begin{equation}\left.\begin{array}{c} \textbf{h}_0(0) = \textbf{a} \\  \textbf{h}_0(t_1) = \textbf{b}  \\ \textbf{h}_1(t_1) = \textbf{b} \\ \textbf{h}_1(t_2) = \textbf{c} \end{array}\right\} \ \mbox{continuity}\\
\end{equation}

That is, we had four unknowns, $\{d_j \}_{j=0}^3$, and used four conditions, all based on continuity, to find them. At this junction, if we wanted to impose more conditions such as continuity across one or more derivatives, we would not have enough degrees of freedom, or free parameters, satisfy all the conditions; we would have an over-constrained system. 

Rather than use linear interpolation, which lead to instantaneous accelerations, let's try to use a cubic polynomial between successive points. Using a higher order polynomial interpolant will also provide more free parameters such that we are able to impose more continuity conditions. Keep in mind, although we will try a cubic polynomial interpolant, our goal is still interpolating between the two states $\textbf{a}=(x_a,y_a)$ and $\textbf{b}=(x_b,y_b).$ 

Our goal is to use a familiar form of an interpolant, that looks awfully reminiscent of the linear case, but with a cubic function of the parameter, $t$, for $t\in[0,1]$. We could attempt to use an interpolant such as the following

\begin{equation}
    \label{cubic_interp} \textbf{h}(t) = \textbf{a} + g(t)( \textbf{b}-\textbf{a}),
\end{equation}
where $g(t)$ is a cubic polynomial, rather than a line as in (\ref{linear_spline}), e.g.,

\begin{equation}
    \label{cubic_interp2} g(t) = d_0 + d_1 t + d_2 t^2 + d_3^3.
\end{equation}
Here we wish for continuity of the function, $\textbf{h}(t)$, continuity in its velocity, $\textbf{h}'(t)$, and no instantaneous accelerations ($\textbf{h}''(t)=0$ at the endpoints of the interpolation domain in $t$). However, when we write the conditions we wish to satisfy,
\begin{align}
    \nonumber
    &\left.\begin{array}{c} \textbf{h}(0) = \textbf{a} \\  \textbf{h}(1) = \textbf{b}  \end{array}\right\} \ \mbox{continuity}\\
    \label{cubic_conds1a}
    &\left.\begin{array}{c} \textbf{h}'(0) = 0 \\  \textbf{h}'(1) = 0  \end{array}\right\} \ \mbox{continuous velocities}\\
    \nonumber
    &\left.\begin{array}{c} \textbf{h}''(0) = 0 \\  \textbf{h}''(1) = 0  \end{array}\right\} \ \mbox{no instantaneous accelerations}
\end{align}
it is clear that we have an over-constrained system, that is, $6$ conditions but only $4$ unknowns, $\{d_j\}_{j=0}^3$. To circumvent this, we can introduce two interpolating mediary points, say $p_1$ and $p_2$, such that we partition the interval $t\in[0,1]$ into three regions: (1) $t\in[0,p_1]$, (2) $t\in[p_1<p_2]$, and (3) $t\in[p_2,1]$. In each of those three regions, we could define an independent cubic interpolant, e.g.,

\begin{equation}
    \label{cubic_interp3} g(t) = \left\{\begin{array}{cc}
        g_0(t)=a_0 + a_1t + a_2t^2 + a_3t^3 \ \ & 0\leq t\leq p_1 \\
        g_1(t)=b_0 + b_1t + b_2t^2 + b_3t^3 \ \ & p_1\leq t \leq p_2 \\
        g_2(t)=c_0 + c_1t + c_2t^2 + c_3t^3 \ \ & p_2 \leq t\leq 1
     \end{array}\right..
\end{equation}

Upon imposing the conditions from (\ref{cubic_conds1a}) onto (\ref{cubic_interp3}), we see that now we have 12 degrees of freedom but only 6 equations, leaving us with an under-constrained system. If we were to think physically about this, at the interfaces $t=p_1$ and $t=p_2$, we would want continuity of our interpolating functions and their first and second derivatives, providing continuity in velocity and acceleration, respectively. Hence the piecewise cubic interpolating functions must satisfy the following constraints:
\begin{align}
    \nonumber
    &\left.\begin{array}{c} g_0(0) = 0 \\  g_2(1) = 1 \\ g_0(p_1)=g_1(p_1) \\ g_1(p_2)=g_2(p_2) \end{array}\right\} \ \mbox{continuity}\\
    \label{cubic_conds2a}
    &\left.\begin{array}{c} g_0'(0) = 0 \\  g_2'(1) = 0 \\ g_0'(p_1)=g_1'(p_1) \\ g_1'(p_2)=g_2'(p_2) \end{array}\right\} \ \mbox{continuous velocities}\\
    \nonumber
    &\left.\begin{array}{c} g_0''(0) = 0 \\  g_2''(1) = 0 \\ g_0''(p_1)=g_1''(p_1) \\ g_1''(p_2)=g_2''(p_2)  \end{array}\right\} \ \mbox{no instantaneous accelerations}
\end{align}

This gives the following linear system to solve, with variables, $p_1$ and $p_2$,

\begin{equation}
    \label{cub_sys} \left[\begin{array}{cccccccccccc} 
    1 &0& 0& 0& 0& 0& 0& 0& 0& 0& 0& 0\\
    0 &1& 0& 0& 0& 0& 0& 0& 0& 0& 0& 0\\
    0 &0& 1& 0& 0& 0& 0& 0& 0& 0& 0& 0\\
    1 & p_1 & p_1^2 & p_1^3 &-1 &-p_1 & -p_1^2 & -p_1^3 & 0 & 0&  0& 0\\
    0 &1 &2p_1& 3p_1^2& 0& -1& -2p_1& -3p_1^2& 0& 0& 0& 0\\
    0 &0 &2 &6p_1& 0& 0& -2& -6p_1& 0& 0& 0& 0\\
    0& 0& 0& 0 &1 &p_2 & p_2^2& p_2^3& -1& -p_2& -p_2^2& -x2^3\\
    0& 0& 0& 0 &0 &1   & 2p_2 & 3p_2^2& 0& -1& -2p_2& -3p_2^2\\
    0& 0& 0& 0& 0 &0   &2 & 6p_2 & 0 & 0 & -2 & -6p_2\\
    0& 0& 0& 0& 0& 0& 0& 0& 1& 1& 1& 1\\
    0& 0& 0& 0& 0& 0& 0& 0& 0& 1& 2& 3\\
    0& 0& 0& 0& 0& 0& 0& 0& 0& 0& 2& 6 \end{array} \right] \left( \begin{array}{c} a_0 \\ a_1 \\a_2 \\ a_3 \\ b_0 \\b_1 \\ b_2 \\ b_3 \\ c_0 \\ c_1 \\ c_2\\ c_3 \end{array}\right) = \left( \begin{array}{c}0\\ 0\\ 0\\ 0\\ 0\\ 0\\ 0\\ 0\\ 0\\ 1\\ 0\\ 0 \end{array}\right) 
\end{equation}

\begin{figure}[H]
        \centering
        \includegraphics[width=0.99\textwidth]{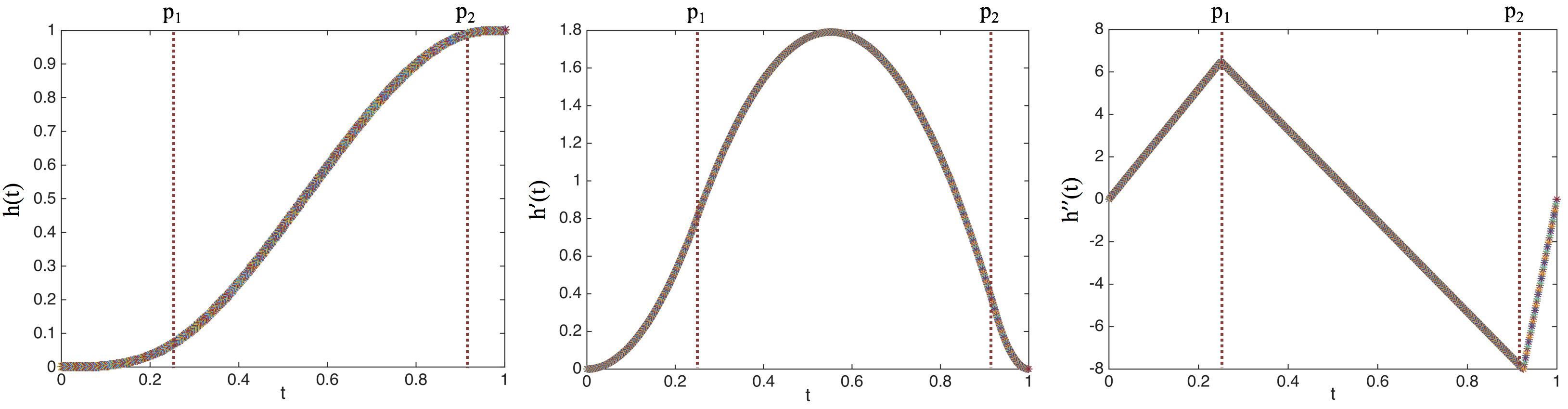}
        \caption{Plots of the piecewise cubic interpolant, $h(t)$, its derivative, $h'(t)$, and its second derivative, $h''(t)$, for $p_1=0.25$ and $p_2=0.925$ for $0\leq t\leq 1$, respectively.}
        \label{Example:poly_f}
    \end{figure}

As an example, if we let $p_1=0.25$ and $p_2=0.925$, upon solving (\ref{cub_sys}), we find the coefficients to be approximately
\begin{equation}
    \label{coeffs} \begin{array}{lll} a_0 = 0 & b_0 = 0.123 & c_0 = -16.778\\ 
    a_1 = 0     & b_1 = -1.481  & c_1 = 53.333 \\
    a_2 = 0     & b_2 = 5.923 & c_2 = -53.333 \\
    a_3 = 4.324\ \ \ \ \ \ & b_3 = -3.577\ \ \ \ \ \ & c_3 = 17.778. \\\end{array}
\end{equation}
A plot of the resulting interpolant, $h(t)$, $h'(t)$, and $h''(t)$ is provided in Figure \ref{Example:poly_f}. It is clear that all the conditions sought after in (\ref{cubic_conds2a}) are satisfied. Moreover by introducing two new parameters $p_1$ and $p_2$, we can essentially control the acceleration of the interpolated motion. The script used to solve this system is provided in the Supplemental Materials, e.g., the \texttt{interp$\_$Function$\_$Coeffs.m} script.
    
As $p_1\rightarrow0$ (or $p_2\rightarrow1$), the initial acceleration (or final deceleration) becomes larger in magnitude. In practice we can use the parameters $p_1$ and $p_2$ to match the acceleration to the kinematics coming from a biological system or engineering system. These parameters $p_1$ and $p_2$ may actually provide a beneficial tool for capturing the correct kinematics of a system in a mathematical model!

Next, in Example \ref{ex:circles2}, we will illustrate qualitative differences in the fluid dynamics when using a cubic interpolant rather than linear interpolant, as is in the previous example. The corresponding source code for this example with a cubic interpolant is found in \texttt{Examples$\_$Education/Interpolation/Moving$\_$Circle/Cubic$\_$Interp}.

%
%

\begin{exmp}
\label{ex:circles2}

In this example we will use the same prescribed motion described in Example \ref{ex:circles}; however, we will use two different interpolation polynomials - one linear and one cubic to interpolate between successive states. Using the cubic interpolant that was determined above, with $p_1=0.25$ and $p_2=0.925$, we ran simulations and compared the results to those when using the linear interpolation scheme.

Simulations were compared at time-points when the circle would be accelerating or decelerating between State $A\rightarrow B$ and the acceleration at the very beginning of State $B\rightarrow C$. This is illustrated in Figure \ref{Example:circles_interp_compare}, where the magnitude of velocity is used to demonstrate qualitative differences in the underlying fluid motion. It is clear that when using different interpolants to prescribe the motion between two states, it can lead to significant differences in the fluid motion. Movies illustrating the dynamical differences are provided in the Supplemental Materials (\texttt{Supplemental/Circles/Linear$\_$Interp or Supplemental/Circles/Cubic$\_$Interp}).

    \begin{figure}[H]
        \centering
        \includegraphics[width=0.99\textwidth]{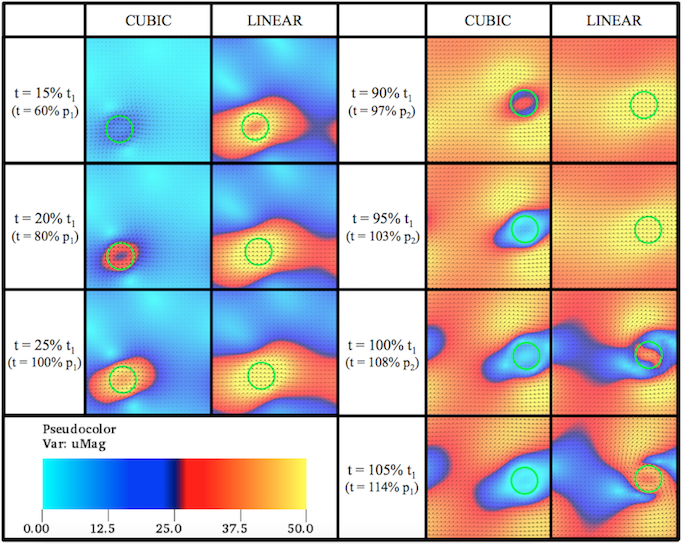}
        \caption{Images illustrating qualitative differences in the magnitude of velocity when using linear and cubic interpolants. Snapshots were taken when the circles were accelerating and decelerating from Phase $A\rightarrow B$, and then accelerating from Phase $B\rightarrow C$.}
        \label{Example:circles_interp_compare}
    \end{figure}

We note that in both cases the circle moves between States $A\leftrightarrow B$ and $B\leftrightarrow C$ with the periods $t_1=0.01$ and $t_2=0.02$, respectively. In fact, qualitatively it appears that in both cases the circles look like they maybe moving in the same way; however, there are clear dynamical differences as seen by the underlying fluid velocity. Again, this is because the velocity and acceleration/deceleration of the circles moving between the states is significantly different. This is an important aspect that should get proper attention when mathematically modeling using prescribed motion. Not only is it important to make the an object begin and end in the right place, but we must also make sure the \textit{way} it moves between the states is biologically (or scientifically) relevant! Introducing higher order polynomial interpolants is a convenient way to introduce more degrees of freedom into a model, so it is able capture more kinematic accuracy.

\end{exmp}


%
%

%
%

\section{Interpolation and beating hearts: a virtual walk through}
\label{Heart_Pumping_Example}
$ $\\

\begin{figure}[H]
        \centering
        \includegraphics[width=0.35\textwidth]{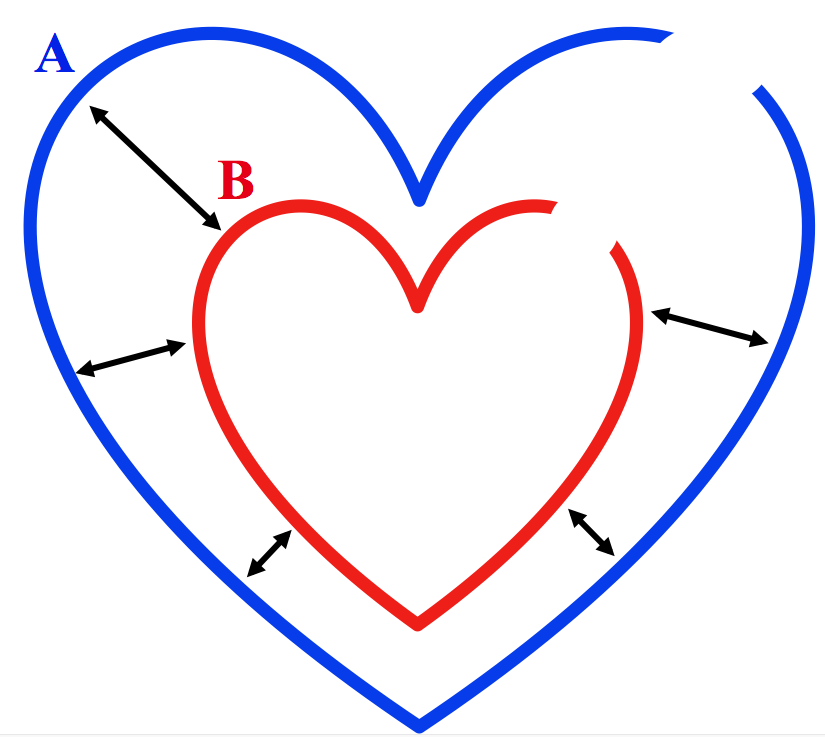}
        \caption{Moving between States A and B to model a beating heart.}
        \label{fig:2StateHearts}
    \end{figure}

Here we present an example of how to implement an object's prescribed motion within the \textit{IB2d} software. We will consider the motion of a beating cartoon heart, that is, a heart that goes between two states, one larger and one smaller, see Figure \ref{fig:2StateHearts}. The hole in the heart is to allow fluid to flow in and out of it,  thereby obeying fluid volume conversation.
    
Running the simulation found in \texttt{Examples$\_$Education/Interpolation/Beating$\_$Heart}, will produce data that can be visualized, as in Figure \ref{fig:BeatingHeart}. The corresponding movie is provided in the Supplemental Materials (\texttt{Supplement/Pulsing$\_$Heart}). We are using the same cubic interpolation scheme that was discussed in Section \ref{Interpolation} to move between State $A\rightarrow B$ and then State $B\rightarrow A$ with periods $t_1$ and $t_2$, respectively. However we also introduce an intermediate resting state of length $t_R$, before moving back from State B$\rightarrow$A to introduce additional possible model complexity.
    \begin{figure}[H]
        \centering
        \includegraphics[width=0.95\textwidth]{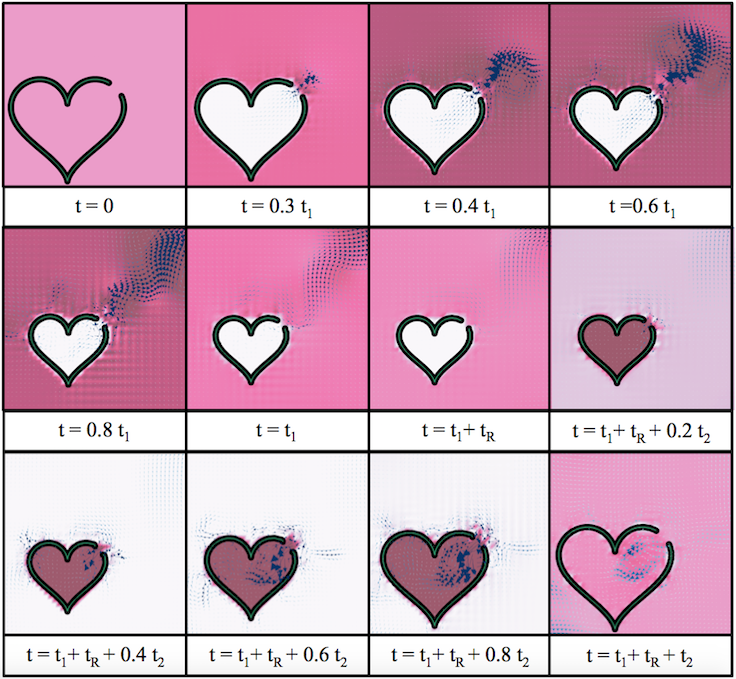}
        \caption{Snapshots of a simulation of a beating cartoon heart that is immersed within a fluid. The colormap depicts the underlying pressure, while the vector field depicts the fluid velocity itself.}
        \label{fig:BeatingHeart}
    \end{figure}
    
We will now dive into detail on how to implement the cubic interpolant to prescribe motion. Although, a beating heart example is introduced here, it should be noted that this will work for just about any geometry, as long as each state has both the same number of points, is ordered consistently, and has a `hole' to obey volume conversation. 

The script that actually prescribes the motion is \texttt{update$\_$Target$\_$Point$\_$Positions.m}. This script does the following three things:
\begin{enumerate}
    \item \textbf{\textit{Specify the period spent moving between states and initialize the cubic interpolant}}.\\
    
        We initialize the time spent in each phase moving between $A\rightarrow B$, resting, and finally $B\rightarrow A$ as $t_1, t_R$, and $t_2$, respectively. We also specify the parameters for the specific cubic interpolant we are going to use to move between States, that is, the coefficients of the cubic interpolant in each sub-phase, $\{a_j,b_j,c_j\}_{j=0}^{3}$, and location of the interpolation nodes, $p_1$ and $p_2$. The values of $p_1$ and $p_2$ were chosen to be $0.25$ and $0.925$, respectively, which is the same case as in Section \ref{Interpolation}.\\ 
        
        Note we also define a period of the total heart beat to be the sum of all the subsequent phases, $t_1+t_R+t_2$, and use modular arithmetic, with respect to said period, for an adjusted time in the simulation in order to simulate repetitive heartbeats.
    
        \begin{figure}[H]
        \centering
        \includegraphics[width=0.9\textwidth]{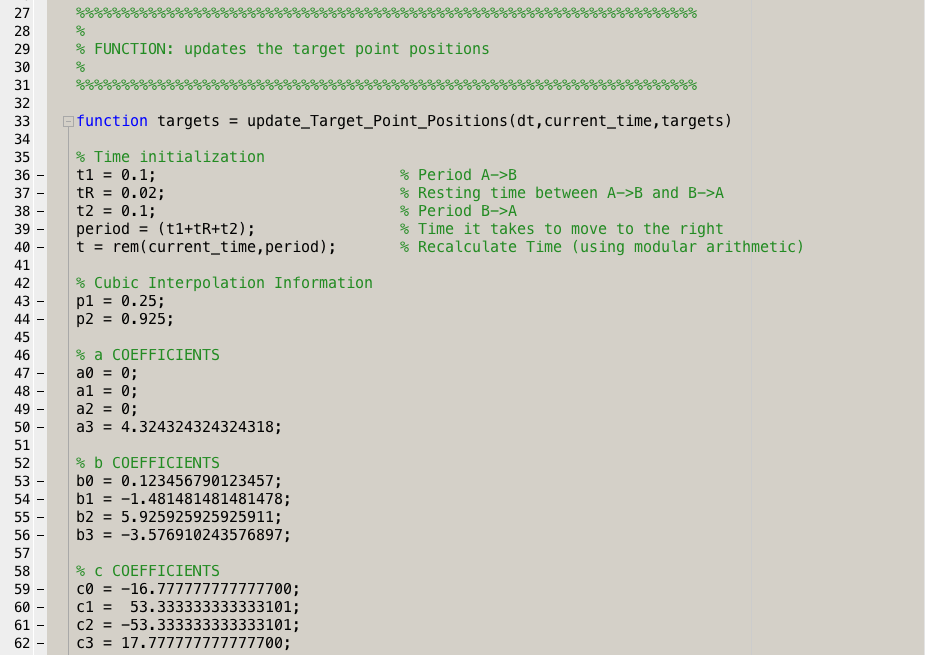}
        \caption{Initializing the time for each phase of motion as well as the cubic interpolant's coefficients from Section \ref{Interpolation}.}
        \label{fig:heart_time_init}
        \end{figure}
        
        $ $\\
    
    \item \textbf{\textit{Read in the points associated for States $A$ and $B$}}.\\
    
        Next we read in the $(x,y)$ positions for each state into $N\times2$-sized matrices, where the columns give the $x$ and $y$ positions, respectively. 
    
        \begin{figure}[H]
        \centering
        \includegraphics[width=0.95\textwidth]{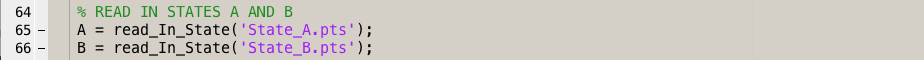}
        \caption{Reading in the (x,y) positions for States $A$ and $B$ into matrices \textbf{A} and \textbf{B}.}
        \label{fig:heart_readAB}
        \end{figure}
        
        For completeness the code that reads in the data from the files \texttt{State$\_$A.pts} and \texttt{State$\_$B.pts} is shown below. 
        
        \begin{figure}[H]
        \centering
        \includegraphics[width=0.95\textwidth]{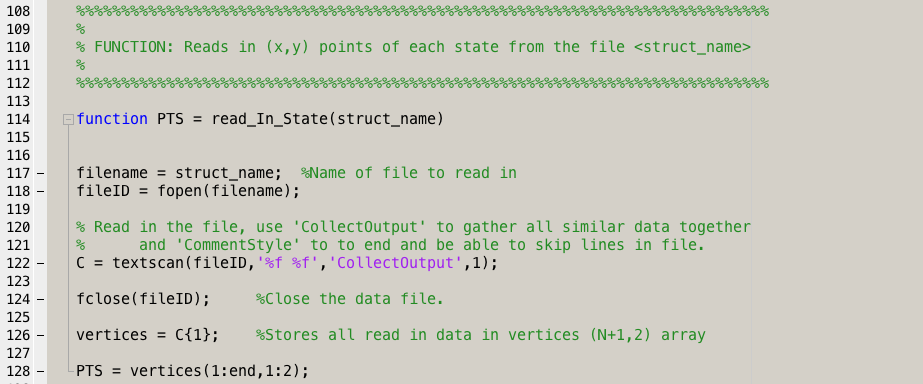}
        \caption{Function that reads in the $(x,y)$ point data.}
        \label{fig:heart_readFunc}
        \end{figure}
        
        We note that the information contained within the files \texttt{State$\_$A.pts} and \texttt{State$\_$B.pts} are lists of the $x$ and $y$ points for each phase, respectively. If you would like to substitute your own shape, rather than use a heart, one only needs to create \textit{.txt} files that contain their own $(x,y)$ point geometries. Note you must also make the \textit{.vertex} file contain the $(x,y)$ positions of the first state as well as include a similarly constructed \textit{.target} file, see the Tutorials in Appendix \ref{app:IB2d_info} for further details.  \\ $ $\\
    
    \item \textbf{\textit{Check which phase of the beating heart it's in, e.g., contraction or expansion, and then update the target point positions to which prescribes the motion of the beating heart}}.\\
    
        Upon checking to see which phase of the simulation the adjusted time currently relates to gives three state possibilities: either the simulation is between States $A\rightarrow B$ or States $B\rightarrow A$, or no motion is being prescribed, e.g., heart is in a rest state.\\
        
        For example, if the simulation time, $t$, is less than the period moving from $A\rightarrow B$, the script then inquiries to find the point between State $A$ and $B$ that it is in, that is, it scales the time appropriately to $\tilde{t} = t/t_1$, so that it is possible to compare $\tilde{t}$ to the interpolation nodes, $p_1$ and $p_2$.
    
        \begin{figure}[H]
        \centering
        \includegraphics[width=0.95\textwidth]{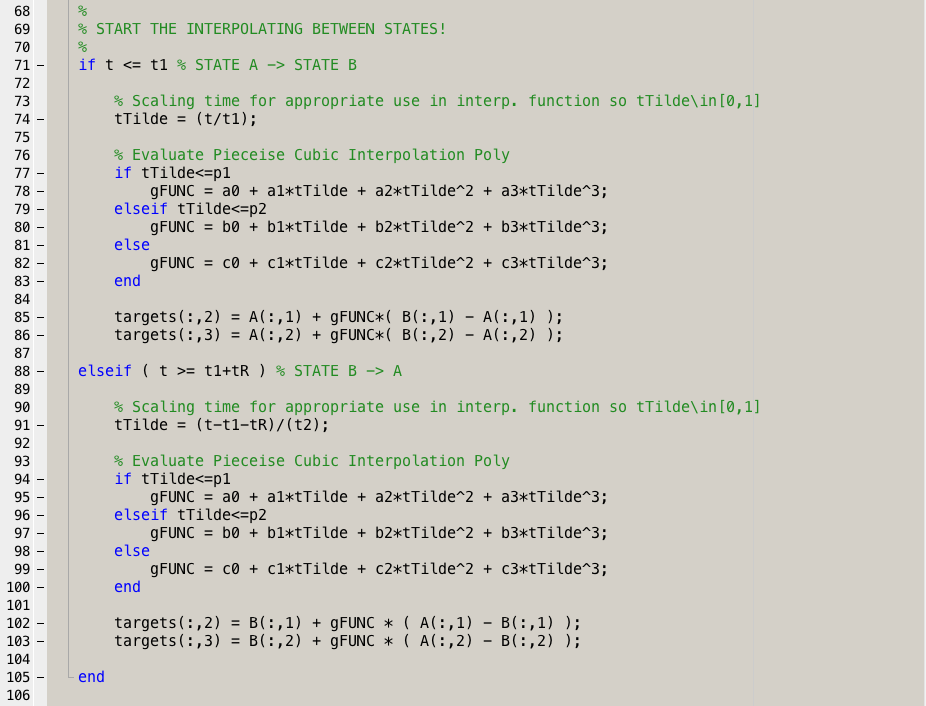}
        \caption{Checks to see which phase of the motion the adjusted simulation time currently relates to and then updates the position of the target points in the $x$ and $y$ directions, which will effectively drive the motion of the beating heart.} 
        \label{fig:heart_code_prescribe}
        \end{figure}
    
\end{enumerate}

%
%
%
%

%
%

\section{Interpolation between material property states: it swims!}
\label{sec:swimmer}

Ready, Set, Swim! Here we present a simple, idealized model of anguilliform locomotion - swimming. Here we do not wish to prescribe the exact kinematics of the swimmer's locomotive patterns, but rather we will only model how the swimmer's body switches between two preferred curvature states. This is a biologically relevant modeling assumption as muscle activation patterns produce specific intrinsic curvatures for a swimmer's body \cite{McMillen:2006,McMillen:2008,Hamlet:2018}. By switching between two different curvature states, the swimmer's body bends and contorts, and locomotion emerges due to the swimmer's interactions with the surrounding fluid. How can model the process of switching between curvature states? That's right; you guessed it - interpolation!

\begin{figure}[H]
    \centering
    \includegraphics[width=0.8\textwidth]{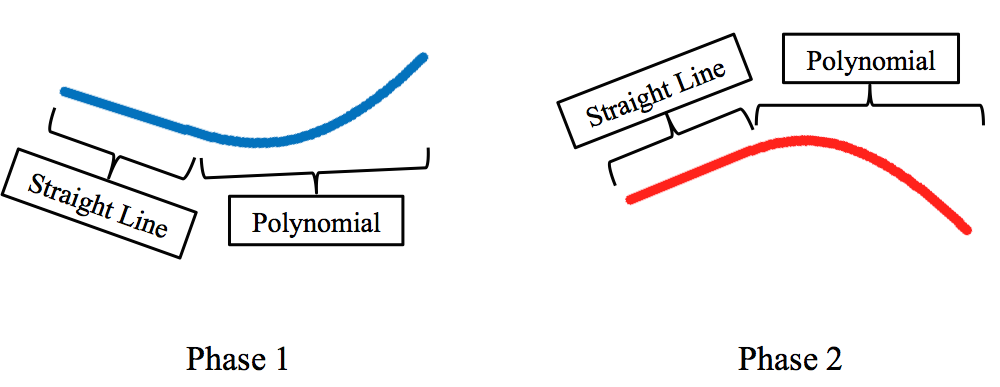}
    \caption{The two phases, in which, the preferred curvature was interpolated between to cause forward swimming, adapted from \cite{BattistaIB2d:2017}.}
    \label{Example:Beam_Swimmer_Phases}
\end{figure}

We must first get in the water before we can swim; let's begin with the shape of the swimmer. To create a simplified scenario, the idealized swimmer's body was constructed by taking a line segment and attaching a polynomial section to it, see Figure \ref{Example:Beam_Swimmer_Phases}, adapted from \cite{BattistaIB2d:2017}. Thus the swimmer's geometry (morphology) is modeled as an infinitely thin $1D$ curve only. The straight portion composes 28\% of the total length of the body, while the polynomial, i.e., $y=x^3$, portion makes up the remaining 72\%. The polynomial section was determined by starting at $x=0$ and adding equally spaced points until $x=L/10$. 

Note that all the points are equally spaced at a distance twice of that of the fluid mesh ($ds=2dx$). Each phase was defined by negating the y-coordinate of the polynomial portion of the body. The ``curvatures" were computed as follows (to tie into the \textit{IB2d} framework, see \cite{BattistaIB2d:2017}):
%
\begin{equation}
    \label{beam:curvature} \begin{array}{c} C_x^P = x_{Lag}^P(s) - 2x_{Lag}^P(s+1) + x_{Lag}^P(s+2) \\ 
    
    C_y^P = y_{Lag}^P(s) - 2y_{Lag}^P(s+1) + y_{Lag}^P(s+2) \end{array}
\end{equation}
where $s$ runs over all Lagrangian points along the swimmer's body and $P$ refers to Phase 1 or 2. 

This intrinsic curvature is the quantity we will now interpolate between. We are no longer interpolating between explicit \textit{positions}, but instead material property states! Although seemingly different, the mathematics (spline interpolation) works out exactly the same. In lieu of changing explicit coordinates (or positions), we now update the curvatures, $C_x^P$ and $C_y^P$ in the \texttt{update$\_$nonInv$\_$Beams.m} script. 

We also define the downstroke and upstroke to be moving between Phase 1 to Phase 2 and Phase 2 to Phase 1, respectively. Furthermore we also define 1 stroke period to encompass both the upstroke and downstroke. The same interpolation rigmarole, as in Section \ref{Heart_Pumping_Example}, follows. 

\begin{figure}[H]
    \centering
    \includegraphics[width=0.995\textwidth]{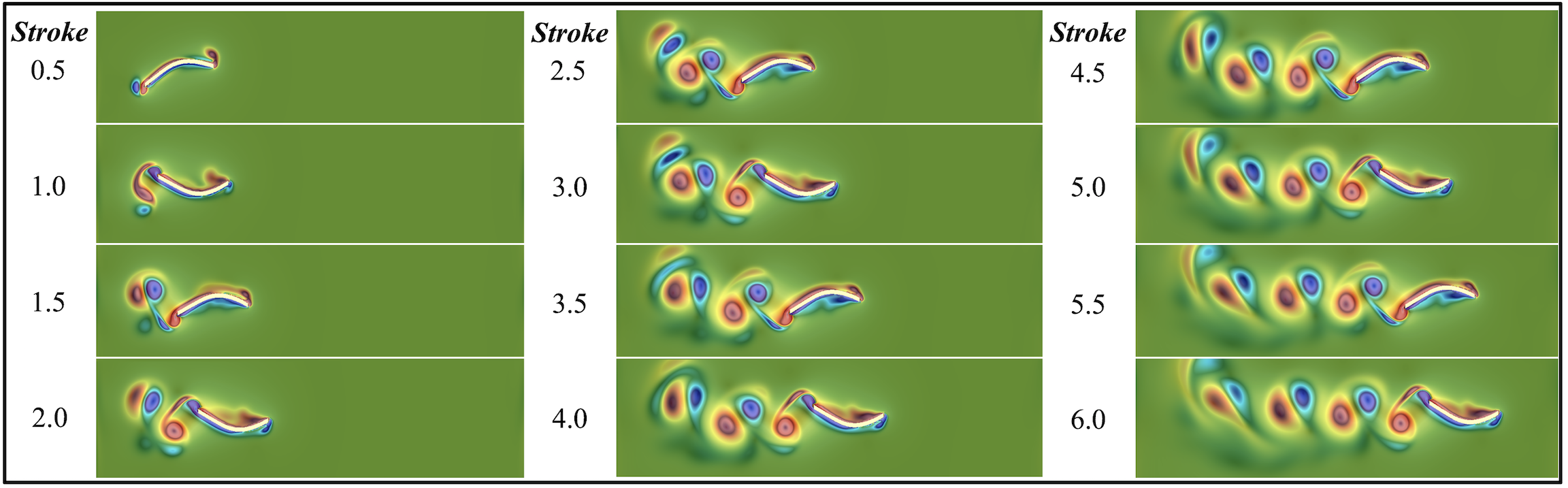}
    \caption{An idealized anguilliform swimmer progressing forward due to continually changes in the preferred curvature of its configuration with a stroke frequency $f=0.5 s^{-1}$. The background colormap illustrates the fluid's vorticity.}
    \label{fig:BeamSwimmer1}
\end{figure}

\begin{figure}[H]
    \centering
    \includegraphics[width=0.875\textwidth]{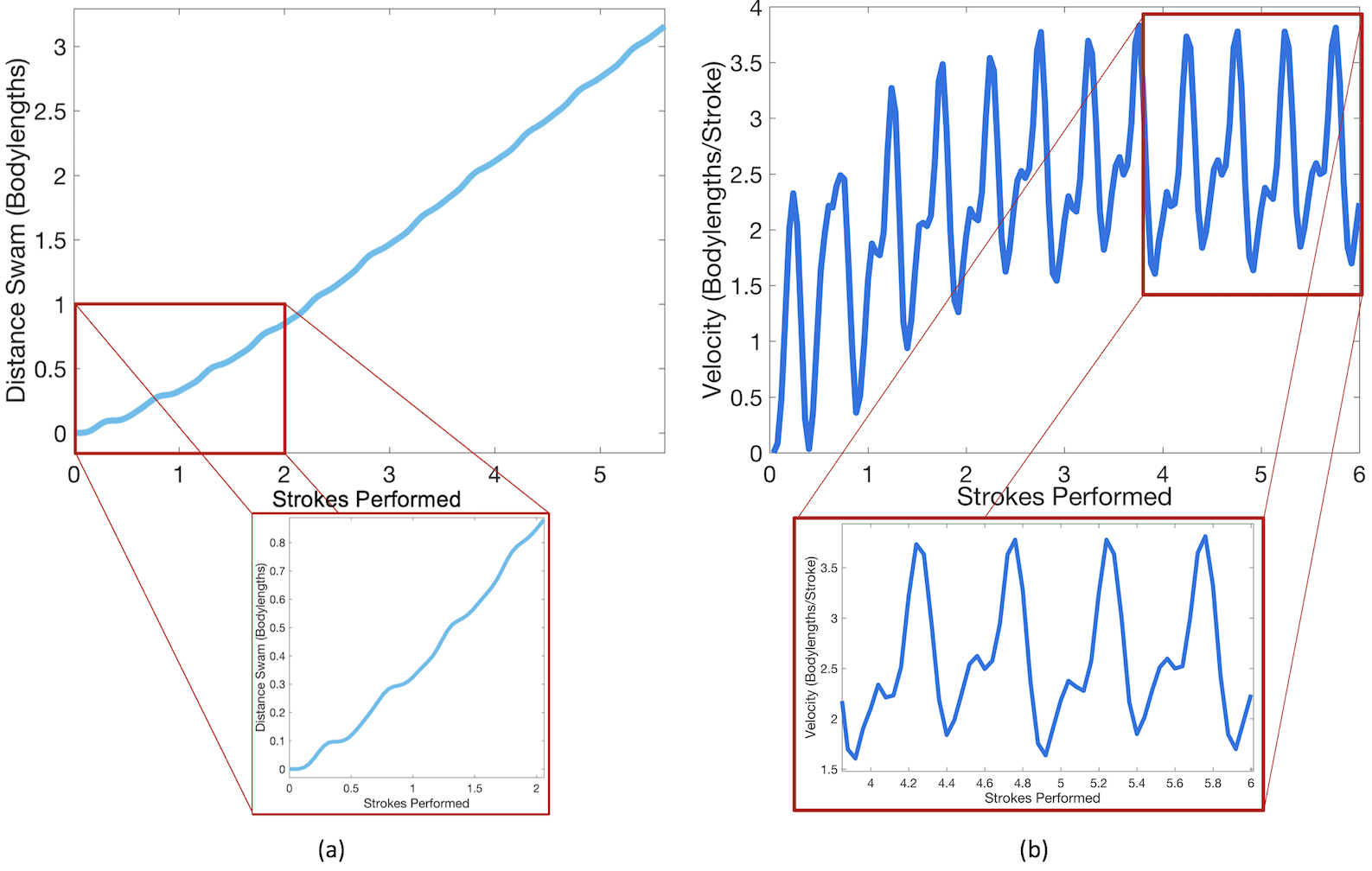}
    \caption{Swimming performance of the single anguilliform swimmer shown in Figure \ref{fig:BeamSwimmer1}. (a) Distance (bodylengths) vs number of strokes performed and (b) velocity (bodylengths/stroke) vs. number of strokes performed.}
    \label{fig:BeamSwimmer1_Distance}
\end{figure}

Running the simulation found in \texttt{Examples$\_$Education/Interpolation/Swimmer/Single$\_$Swimmer} will produce locomotion data that can be visualized as in Figure \ref{fig:BeamSwimmer1}. This figure shows the idealized anguilliform swimmer moving forward due to vortices being shed off its caudal end during each stroke. The background colormap represents the fluid's vorticity, e.g, the local swirling motion of the fluid (mathematically given by the curl of the velocity field, $\nabla\times\textbf{u}({\bf x},t)$). The corresponding movie to Figure \ref{fig:BeamSwimmer1} is provided in the Supplementary Materials (\texttt{Supplemental/Swimmer/Individual$\_$Swimmer/}). Furthermore, we can quantitatively track the position of the swimmer's head over time, using the script  \texttt{Individual$\_$Swimmer$\_$Analysis.m}, to see what its forward swimming patterns (and performance) looks like, see Figure \ref{fig:BeamSwimmer1_Distance}. 

At this point while we have a single simulation of one anguilliform swimmer, there are many interesting questions one could ask, including a plethora of interesting biological questions. However, we will first focus on how subtle changes in interpolating between curvature states affects swimming performance. Note that for the simulation shown in Figure \ref{fig:BeamSwimmer1} that $(p_1,p_2)=(0.1,0.9).$ In particular, we will ask three questions:
\begin{enumerate}
    \item What happens when the interpolation mediary points $p_1$ and $p_2$ are changed? Remember these points help dictate the acceleration and velocity profile of the interpolation (see Section \ref{Interpolation}).
    \item What happens if we make the interpolation mediary points $(p_1,p_2)$ asymmetric (e.g., say if $p_1=0.1$ and $p_2=0.5$ rather than $p_2=0.9$)?
    \item What happens if we have an asymmetric stroke pattern? (For example, if the upstroke is 25\% of the total period while the downstroke is only 75\%?)
\end{enumerate}

Lastly, we can have a little fun with our swimmer, taking advantage of the fact it is immersed in a fluid, and ask how does changing the fluid environment affect swimming performance? To change its fluid environment, we will only have to vary the fluid's viscosity. This effectively asks how the swimmer performs in stickier and stickier fluid environments, like going from water to corn syrup. For those with previous experience in fluid dynamics, this equates to looking how swimming performance varies over a range of Reynolds Numbers, $Re$.

It is important to note that while asking these questions (and hopefully making hypothesis) we are only changing one parameter of a single simulation at a time, whether that it is $(p_1,p_2)$, the upstroke and downstroke percentages of the total period, or the fluid's viscosity.

%
%

\subsection{Changing $(p_1,p_2)$ symmetrically}
\label{sec:swim_case1}

First we will investigate how the choice of interpolation mediary points $(p_1,p_2)$ affects swimming performance of our idealized anguilliform swimmer. These simulations are found in \texttt{Examples$\_$Education/Interpolation/Swimmer/Case$1$}. We will vary the $(p_1,p_2)$ points symmetrically about the interpolation interval and consider the following cases:
\begin{enumerate}
    \item $(p_1,p_2)=(0.1,0.9)$
    \item $(p_1,p_2)=(0.2,0.8)$
    \item $(p_1,p_2)=(0.3,0.7)$
    \item $(p_1,p_2)=(0.4,0.6)$
\end{enumerate}

Upon varying these points, we need to make sure that our interpolation function is consistent, that is, we need to solve the linear system described in Section \ref{Interpolation} accordingly to get the proper coefficients for the spline interpolant. These coefficients are listed in Supplement 2 of the Supplementary Materials. Once calculated, we can modify the \texttt{update$\_$nonInv$\_$Beams.m} script, which performs the curvature interpolation. 

\begin{figure}[H]
    \centering
    \includegraphics[width=0.99\textwidth]{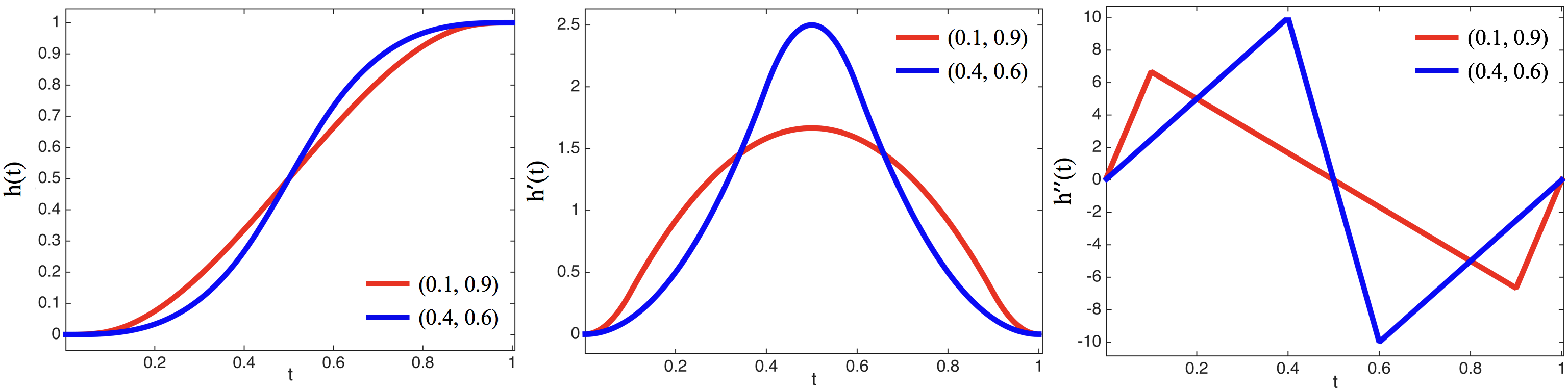}
    \caption{Plots of the piecewise cubic interpolant, $h(t)$, its derivative, $h'(t)$, and its second derivative, $h''(t)$, with $0\leq t\leq 1$, for varying $(p_1,p_2)$ symmetrically chosen.}
    \label{fig:SwimmerSymP1P2_g}
\end{figure}

We will now compare the interpolation profiles ($h(x)$, $h'(x)$, and $h''(x)$) for two cases: $(p_1,p_2)=\left\{(0.1,0.9),(0.4,0.6)\right\}.$ Comparison plots are given in Figure \ref{fig:SwimmerSymP1P2_g}. We note that in every case we still have continuous first and second derivatives; however, the velocity and acceleration profiles are significantly different.

\begin{figure}[H]
    \centering
    \includegraphics[width=0.99\textwidth]{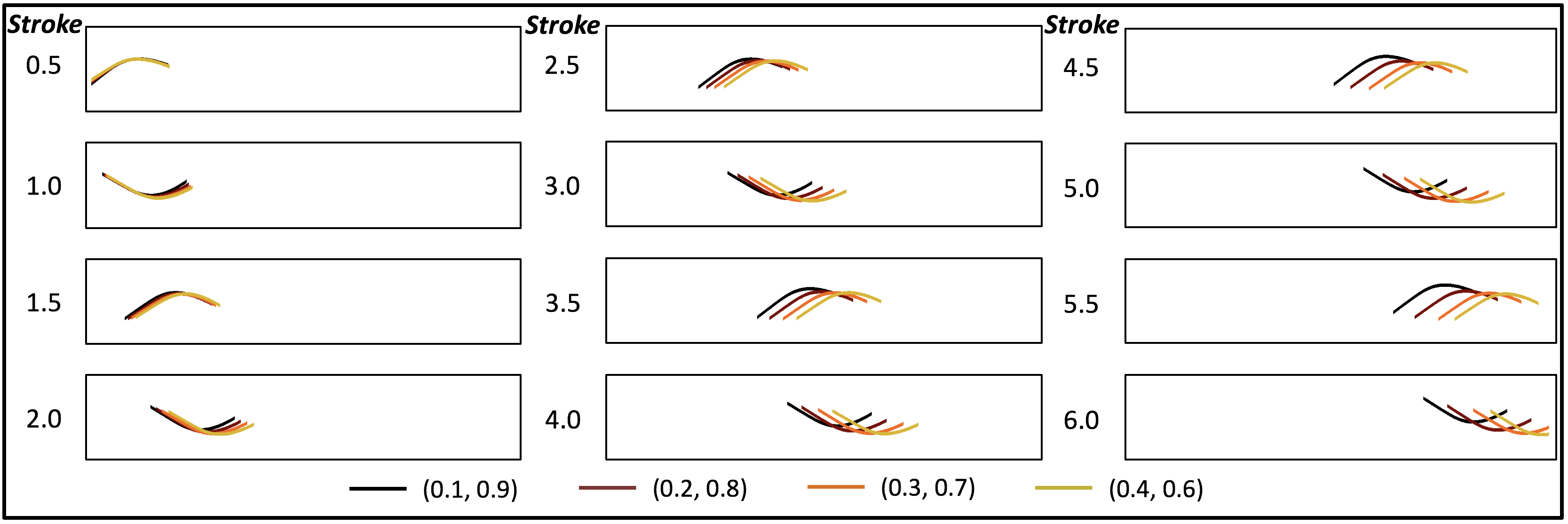}
    \caption{Snapshots from simulations for the case of symmetric interpolation points, given by $(p_1,p_2)\in\{(0.1,0.9),(0.2,0.8),(0.3,0.7),(0.4,0.6)\}$.}
    \label{fig:SwimmerSymP1P2}
\end{figure}

\begin{figure}[H]
    \centering
    \includegraphics[width=0.99\textwidth]{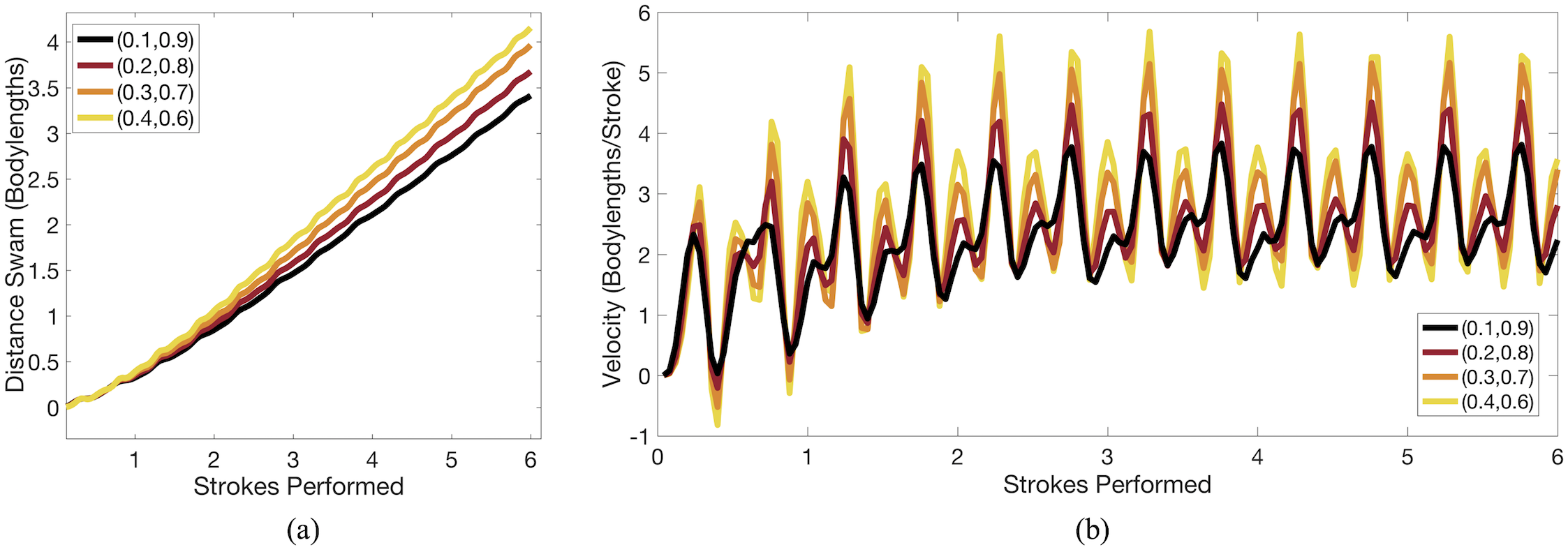}
    \caption{(a) Forward distance swam and (b) forward velocity vs. swimming strokes performed in the case of symmetric interpolation points $(p_1,p_2)$ in $[0,1]$.}
    \label{fig:SwimmerSymP1P2_data}
\end{figure}

Upon running the aforementioned simulations, it is evident that changing ($p_1,p_2$) affects swimming performance! Snapshots from the simulation are given in Figure \ref{fig:SwimmerSymP1P2}. Note that although the swimmer's position from each case are over laid on each other, each simulation was independently performed; there are no swimmer-swimmer interactions. The case when $(p_1,p_2)=(0.4,0.6)$ appears in the lead after $6$ strokes followed by cases $(0.3,0.7), (0.2,0.8)$, and then $(0.1,0.9)$, respectively. The faster cases correspond to higher magnitudes of velocity and acceleration, see Figure \ref{fig:SwimmerSymP1P2_g}. We also present the distance swam vs. swimming stroke as well as forward swimming speed vs. stroke in Figure \ref{fig:SwimmerSymP1P2_data}, which further confirms those results. Furthermore, both peaks in the forward swimming speed's waveform are also higher in the faster cases. The corresponding movie for these simulations is provided in the Supplementary Materials (\texttt{Supplemental/Swimmer/Case1/}).

Simply changing the interpolation mediary points, $p_1$ and $p_2$, affects swimming performance even when everything else remains the same - the same cubic spline-based interpolating function , the same upstroke and downstroke periods, and the same fluid environment! Next we will once again ask how swimming performance is affected if we again change the interpolation points $p_1$ and $p_2$, but this time place them asymmetrically about the interpolation window $[0,1]$.

%
%

\subsection{Changing $(p_1,p_2)$ asymmetrically}
\label{sec:swim_case2}

Here we will again will inquire into how changing the interpolation mediary points ($p_1,p_2$) affects swimming performance, but this time choose $p_2$ such that interpolation points are not symmetric within the interpolation interval $[0,1]$. These simulations are found in \texttt{Examples$\_$Education/Interpolation/Swimmer/Case2}. We selected the following $(p_1,p_2)$ cases:
\begin{enumerate}
    \item $(p_1,p_2)=(0.1,0.9)$
    \item $(p_1,p_2)=(0.1,0.7)$
    \item $(p_1,p_2)=(0.1,0.5)$
    \item $(p_1,p_2)=(0.1,0.3)$
\end{enumerate}

It is important to note that in this section, although we are asymmetrically varying $p_2$ about the interpolation interval, both the upstroke and downstroke have the same period. The only difference is that the rate of change of the interpolating function $h(t)$ during each portion of the stroke.

Again, to ensure that the interpolation function is consistent, we solve the linear system described in Section \ref{Interpolation} for each different set of interpolation points, $(p_1,p_2)$. These coefficients are listed in Supplement 2 of the Supplementary Materials and are used in each corresponding \texttt{update$\_$nonInv$\_$Beams.m} script to perform the curvature interpolation. 

\begin{figure}[H]
    \centering
    \includegraphics[width=0.99\textwidth]{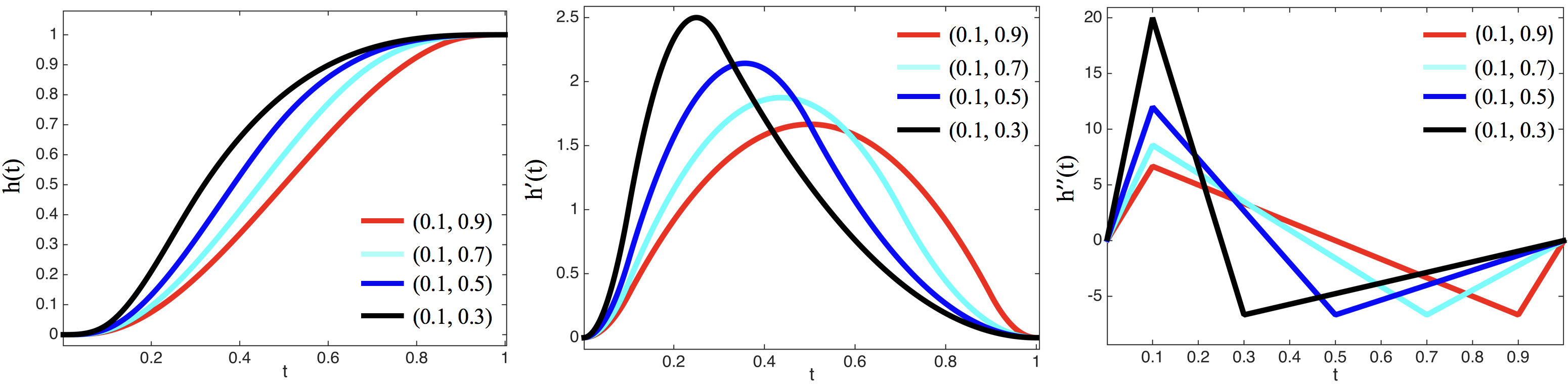}
    \caption{Plots of the piecewise cubic interpolant, $h(t)$, its derivative, $h'(t)$, and its second derivative, $h''(t)$, with $0\leq t\leq 1$, for varying $(p_1,p_2)$ asymmetrically chosen.}
    \label{fig:SwimmerAsymP1P2_g}
\end{figure}

The interpolation profiles $h(t)$, $h'(t)$, and $h''(t)$ look strikingly different than those shown in Section \ref{sec:swim_case1} due to the asymmetry introduced by choice of $p_1$ and $p_2$. The profiles are given in Figure \ref{fig:SwimmerAsymP1P2_g}. 

\begin{figure}[H]
    \centering
    \includegraphics[width=0.99\textwidth]{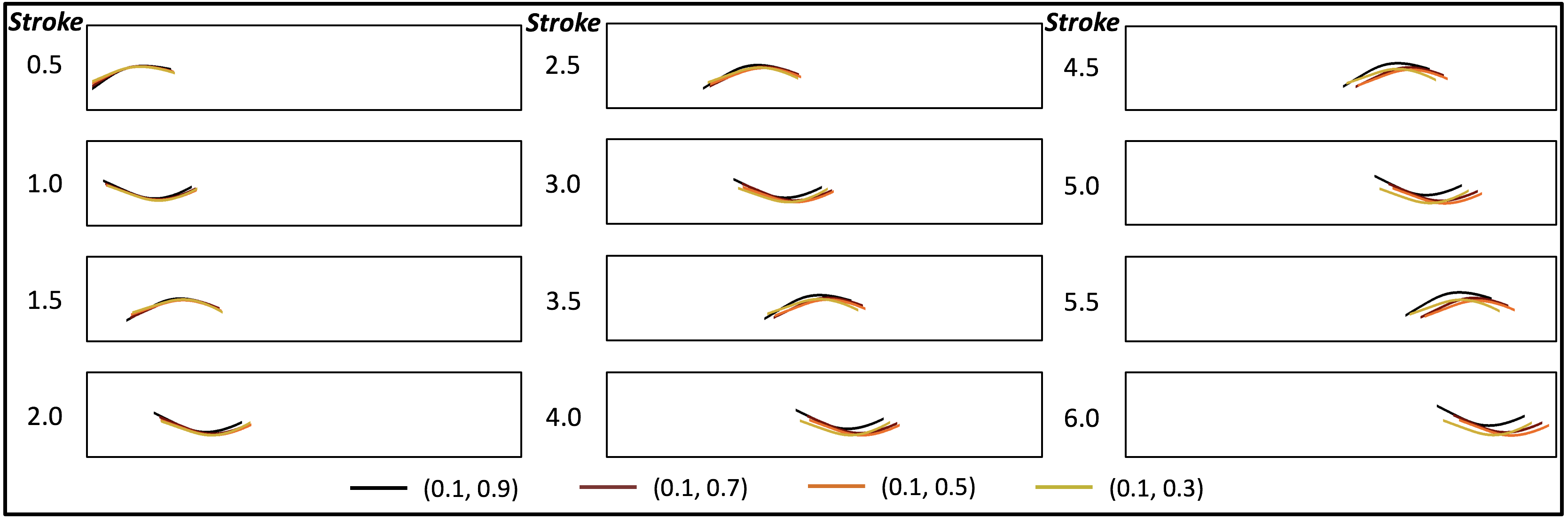}
    \caption{Snapshots from simulations for the case of asymmetric interpolation points, given by $p_1=0.1$ and $p_2\in\{0.3,0.5,0.7,0.9\}$.}
    \label{fig:SwimmerAsymP1P2}
\end{figure}

\begin{figure}[H]
    \centering
    \includegraphics[width=0.99\textwidth]{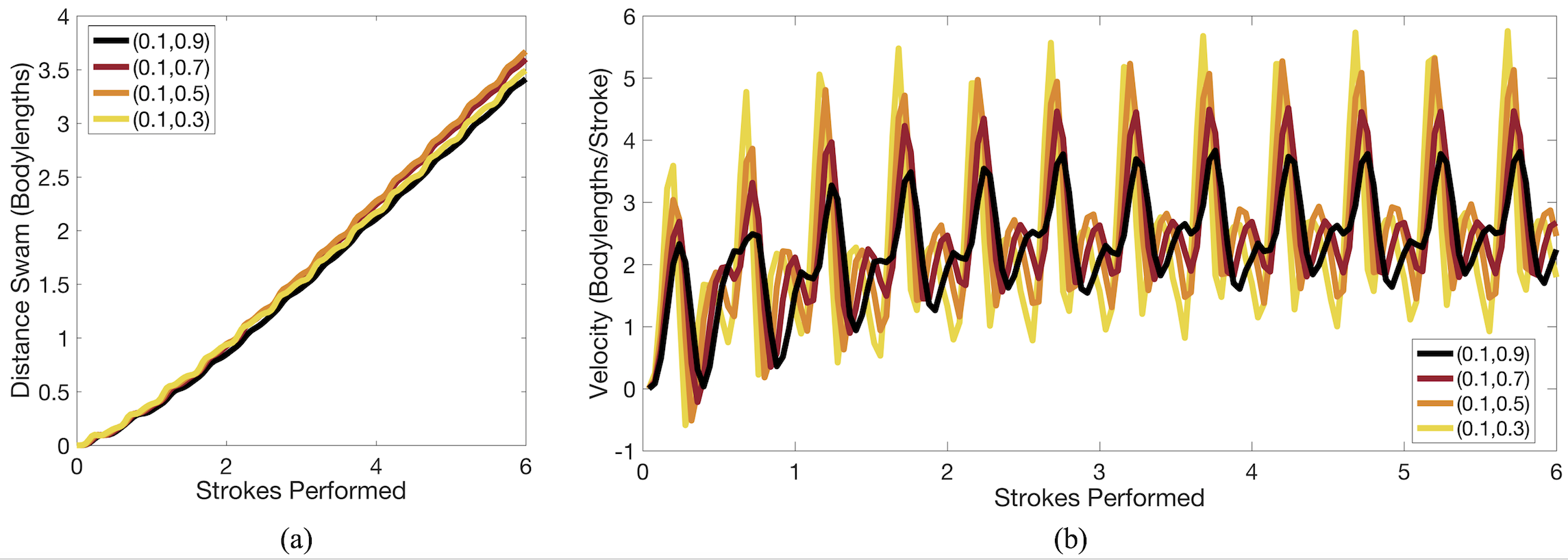}
    \caption{(a) Forward distance swam and (b) forward velocity vs. swimming strokes performed in the case of asymmetric interpolation points, given by $p_1=0.1$ and $p_2=\{0.3,0.5,0.7,0.9\}$.}
    \label{fig:SwimmerAsymP1P2_data}
\end{figure}

As hopefully hypothesized, the dynamics are different between each swimmer for the above cases; however, perhaps surprisingly, there appears to be less variation than the previous case of symmetric $(p_1,p_2)$ choices in terms of forward swimming performance. Snapshots of the four swimmers are shown in Figure \ref{fig:SwimmerAsymP1P2}. In this case there was a non-linear relationship with choice of $p_2$ and how fast the swimmer went, e.g., the case with $p_2=0.5$ was the fastest, followed by $p_2=0.7$, then $0.3$, and finally $0.9$. This is confirmed when analyzing the data, shown in Figure \ref{fig:SwimmerAsymP1P2_data}, which gives the distance swam vs. swimming stroke as well as forward swimming velocity vs. stroke. The corresponding movie of these simulations is provided in the Supplementary Materials (\texttt{Supplemental/Swimmer/Case2/}). What do you think happens if we again sweep over $p_2=\{0.3,0.5,0.7,0.9\}$ but choose a different $p_1$, where $p_1\in(0,p_2)$?

While Sections \ref{sec:swim_case1} and \ref{sec:swim_case2} used different interpolation mediary points, $p_1$ and $p_2$, they both used the same upstroke and downstroke periods as well as same fluid environment, e.g., fluid density and viscosity were the same. We will now investigate variances in swimming performance due to varying stroke periods, followed by changing the fluid environment via varying the fluid's viscosity.

%
%

\subsection{Making asymmetric stroke periods}
\label{sec:swim_case3}

In this case we will keep the interpolation points fixed at $(p_1,p_2)=(0.1,0.9)$ and fix the stroke period to $T=2.0s$ (frequency of $0.5$ Hz). We then asymmetrically vary the upstroke (UPS) and downstroke (DWS) percentages of the total stroke period ($T$). Recall that earlier we defined one stroke to be the upstroke and downstroke periods added together. To that end, we simulated the following cases:
\begin{enumerate}
    \item UPS = DWS, e.g., (UPS,DWS)=(50\%T,50\%T)
    \item UPS = 75\% DWS, e.g., (UPS,DWS)=(42.9\%T,57.1\%T)
    \item UPS = 50\% DWS, e.g., (UPS,DWS)=(33\%T,0.66\%T)   
    \item UPS = 25\% DWS, e.g., (UPS,DWS)=(20\%T,0.80\%T)
\end{enumerate}

Note that although we have made each portion of a single full stroke have a different sub-period, we can still use the same piecewise interpolant, $h(t)$, to interpolate between each! These simulations are found in \textit{Examples$\_$Education/Interpolation/Swimmer/Case3}.

\begin{figure}[H]
    \centering
    \includegraphics[width=0.99\textwidth]{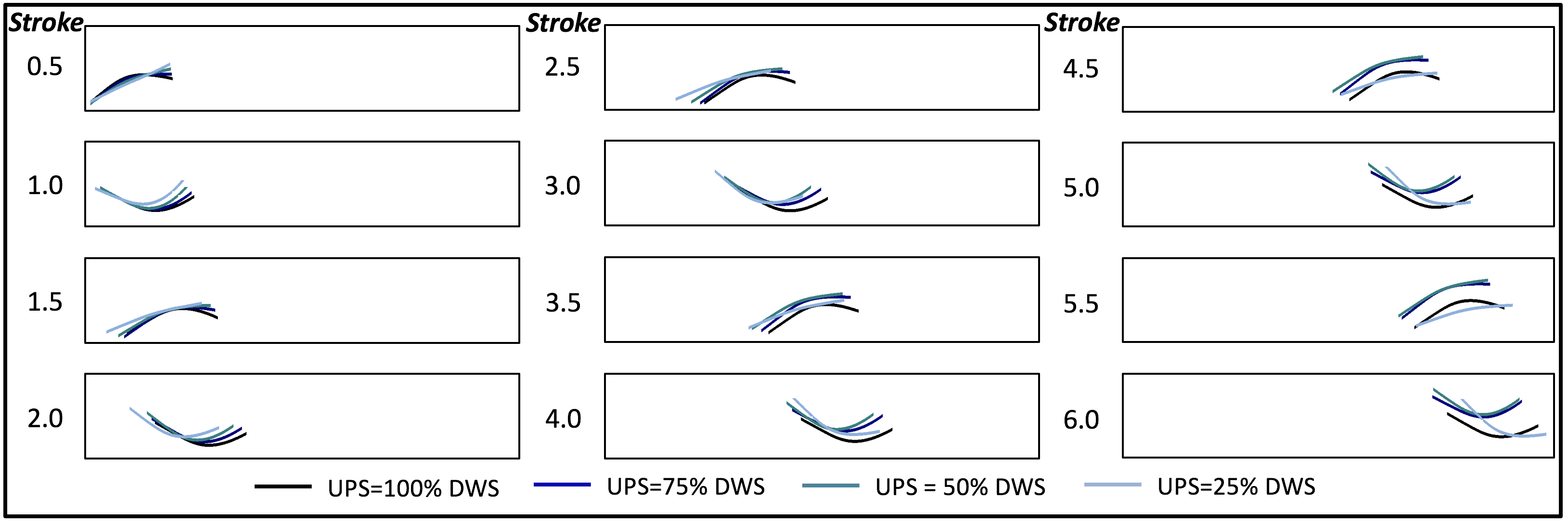}
    \caption{Snapshots from simulations with varying upstroke and downstroke percentages of a single stroke period.}
    \label{fig:SwimmerAsymStrokes}
\end{figure}

\begin{figure}[H]
    \centering
    \includegraphics[width=0.99\textwidth]{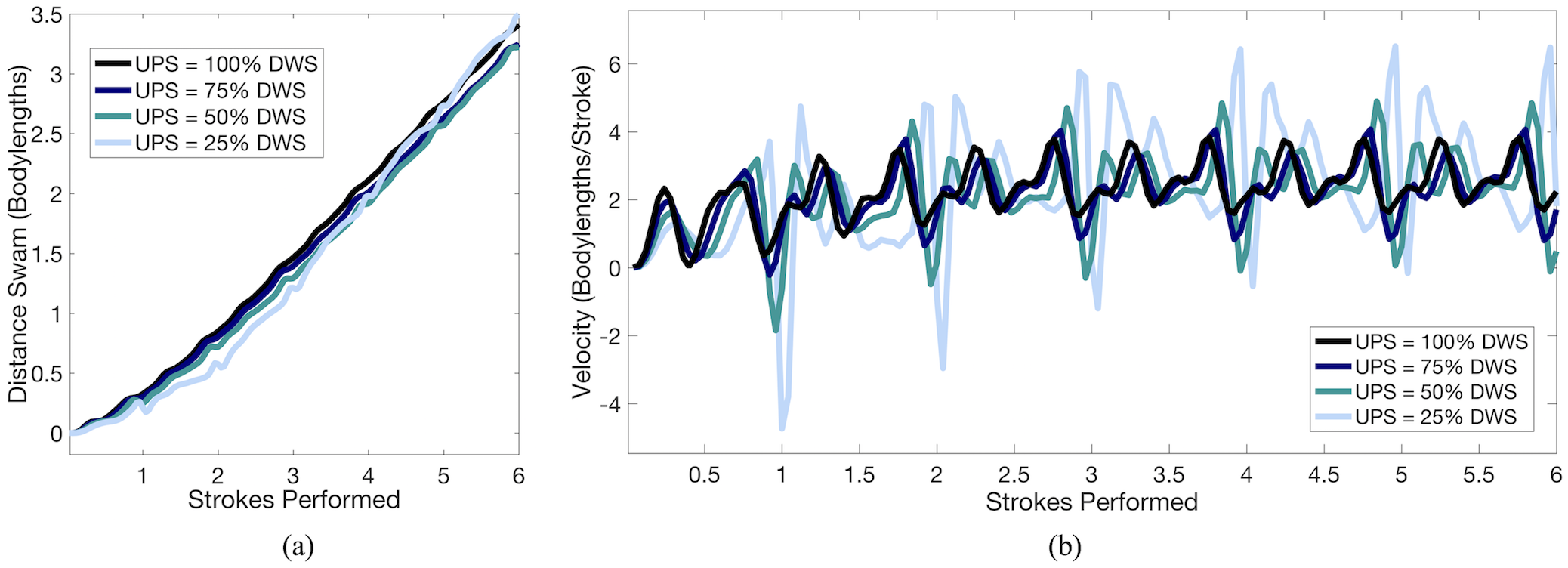}
    \caption{(a) Forward distance swam and (b) forward velocity vs. swimming strokes performed in the case of asymmetric upstroke and downstroke periods.}
    \label{fig:SwimmerAsymStrokes_data}
\end{figure}

As the UPS percentage of a stroke decreases, the upstroke happens faster. However, although the swimmer that swims forward the fastest also has the quickest UPS, having a faster UPS does not always equate to a faster forward swimming speed, see Figures \ref{fig:SwimmerAsymStrokes} and \ref{fig:SwimmerAsymStrokes_data}. The initial acceleration of the UPS=25\%DWS case is the slowest but eventually it starts outswimming the others - truly a tortoise and a hare story (well not exactly, biologically). Figure \ref{fig:SwimmerAsymStrokes} gives snapshots of the four swimmers and Figure \ref{fig:SwimmerAsymStrokes_data} presents the distance swam vs. swimming stroke as well as forward swimming velocity vs. swimming stroke. The corresponding movie of these simulations is provided in the Supplementary Materials  (\texttt{Supplemental/Swimmer/Case3/}). Interestingly, due to the asymmetric UPS and DWS, the swimming velocity profiles look significantly different than those in Figures \ref{fig:SwimmerSymP1P2_data} and \ref{fig:SwimmerAsymP1P2_data}. In particular, the waveforms appear trimodal rather than bimodal, which were observed in the cases of varying $(p_1,p_2)$, especially in the cases of UPS = 25\% DWS and UPS = 50\% DWS. 

What do you think would happen if we redid this same analysis, but with a different $(p_1,p_2)$? Or if we varied the stroke frequency cycle-by-cycle during the simulation?

%
%

\subsection{Changing the fluid viscosity ($Re$)}

Finally, we will consider what happens if we put the swimmer in varying fluid environments, via changing the fluid's viscosity. This equates to placing the swimmer in less or more of a viscous fluid. Examples of highly viscous fluids include things like honey or corn syrup, or fluids that are generally ``thicker" or "more sticky", while less viscous fluids, like water, are considerably less so. For these numerical experiments we keep all other parameters the same, i.e., all the interpolation parameters, upstroke and downstroke periods, geometry, etc. We considered fluid dynamic viscosities, $\mu$, across $5$ orders of magnitude from $0.05$ to $5000$. Note that the viscosity considered in all previous cases (Sections \ref{sec:swim_case1}-\ref{sec:swim_case3}) was $\mu=10.$

As briefly stated earlier, this is equivalent to varying the Reynolds Number, $Re$, which describes the ratio of inertial to viscous forces, which is quantitatively given by
\begin{equation}
    \label{eq:Re} Re = \frac{\rho VL}{\mu}.
\end{equation}
Note that $\rho$ and $\mu$ are the fluid's density and dynamic viscosity, respectively, while $L$ and $V$ are characteristic length and velocity scales for the system. We will not go into more depth regarding Reynolds Number; more information regarding $Re$ ``scaling" studies can be found in \cite{Borazjani:2008,Hershlag:2011,Baumgart:2014,Battista:2016a,Battista:2019,Miles:2019b}. Let's see how these idealized swimmers perform in different viscosities!

\begin{figure}[H]
    \centering
    \includegraphics[width=0.99\textwidth]{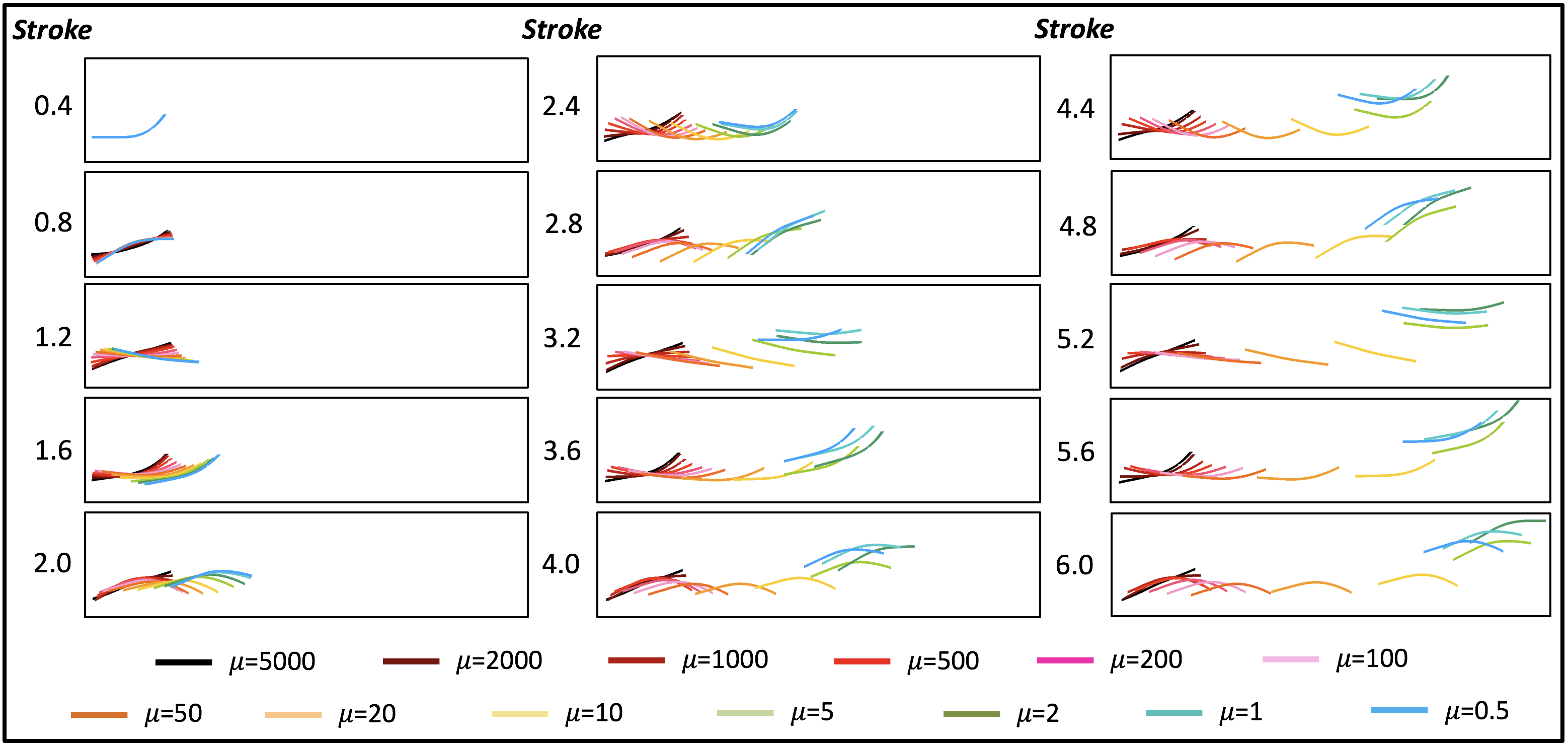}
    \caption{Snapshots from simulations with varying fluid viscosities.}
    \label{fig:SwimmerViscosity}
\end{figure}

\begin{figure}[H]
    \centering
    \includegraphics[width=0.99\textwidth]{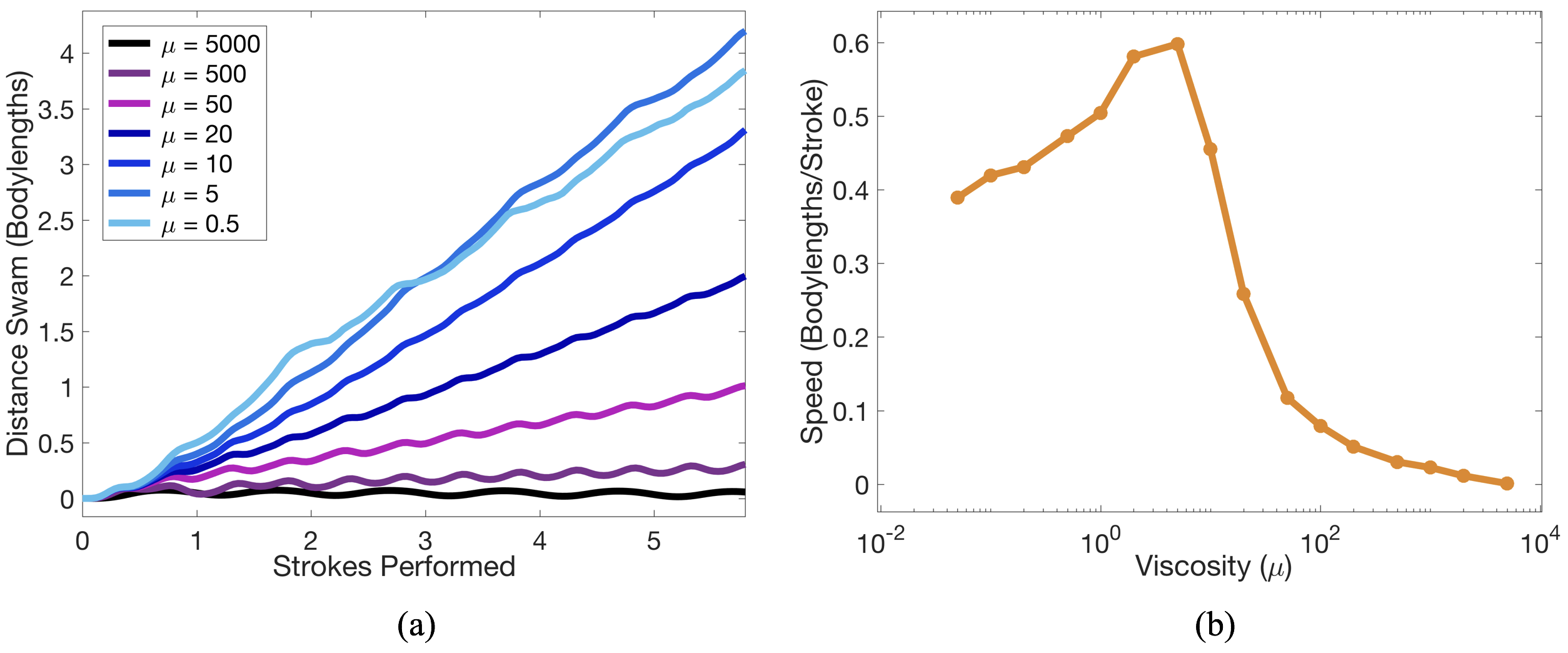}
    \caption{(a) Forward distance swam vs swimming strokes performed and (b) swimming speed (bodylengths/stroke) vs. viscosity.}
    \label{fig:SwimmerViscosity_data}
\end{figure}

Snapshots from simulations of various swimmers in fluids with different viscosities are provided in Figure \ref{fig:SwimmerViscosity}. The corresponding movie is provided in the Supplementary Materials (\texttt{Supplement/Swimmer/Viscosity$\_$Race/}). Qualitatively it appears that swimming performance of our idealized anguilliform swimmer decreases as viscosity increases. When the fluid is ``thick" or ``sticky"-enough, the swimmer may not even able to move forward with this set of model parameters (see the $\mu=5000$ case) unlike its anguilliform counterparts in less viscous fluid! This is confirmed in Figure \ref{fig:SwimmerViscosity_data}, which gives the distance swam (bodylengths) vs. swimming strokes performed and average forward swimming speed (bodylengths/stroke) vs viscosity ($\mu$). Interestingly, it appears that this particular anguilliform swimmer has a maximum speed at a particular viscosity around $\mu\sim500$. That is, in this model of anguilliform locomotion, simply putting the swimmer into less and less viscous fluid will not always result in a faster swimming speed. How do you think this would change if you varied some of the interpolation parameters, $(p_1,p_2)$, or the stroke frequency?

\section{Discussion}
\label{discussion}

Hopefully this has convinced you that there are some practical uses of interpolation in mathematical modeling, which are not generally discussed in traditional numerical analysis settings. In this paper we illustrated a few of the possibilities when applying spline interpolation techniques to mathematical modeling, including prescribing movement patterns (Sections \ref{Interpolation} and \ref{Heart_Pumping_Example}) and material property states (Section \ref{sec:swimmer}). In particular, we demonstrated the following practical aspects of interpolation in mathematical modeling:
\begin{enumerate}
    \item Interpolation can be used to prescribe the motion of an object.
    \item Interpolation can be used to switch between different material property states of an object, which can give rise to unsuspecting, interesting dynamics.
    \item When using spline interpolants, the number of continuous derivatives affects the resulting dynamics of the system. That is, it does not only matter that you get from \textit{A} to \textit{B}, but also \textit{how you get there}, in terms of velocities and accelerations. 
    \item Thus to relinquish modeling artifacts, one could design their interpolant to match observed velocities and accelerations from experimental data, if possible.
    \item Even when not prescribing the precise movement of an object, but rather the object's material property states (e.g., curvature), changing the spline interpolant affects the system's outcome. 
    \item In fact, subtly changing aspects of the interpolant can lead to significant changes in the unveiling dynamics.
\end{enumerate}

We note that the simulations in Sections \ref{Interpolation} and \ref{Heart_Pumping_Example} were designed on a coarse mesh so that students can run them locally on laptops in a manner of a few minutes. However the swimmer simulations in Section \ref{sec:swimmer} were constructed on much finer meshes, which have been observed to be required for locomotion previously \cite{BattistaMizuhara:2019}. Each of the swimmer simulations takes on the order of $\sim2$ hours on a personal machine ($\sim$4-16GB RAM, $\sim$2-3GHz processor). In all of these examples, students have the opportunity to experience scientific computing research in practice, e.g., simulations that can greatly vary in computational time, produce a lot of data with non-trivial data analysis, and open the floor for discussions on effective data visualization. 

The main purpose of this work was to bring interpolation to life for students, allowing them to visually witness how subtle differences in interpolation techniques can lead to significant differences in dynamics, particular within mathematical models. For this reason all codes, both simulation and analysis scripts, are made available. To that extent, this work allows students the opportunity to ask a variety of questions (e.g., such as those posed in Section \ref{sec:swimmer}), explore, and chase their answers. This encourages students to `play' in a numerical and mathematical setting, experiencing mathematical material in a possibly unfamiliar way. Francis Su, former MAA President, has publicly said, ``\textit{Play is part of  human flourishing. You cannot flourish without play. And if mathematics is for human flourishing, we should ``play up" the role of play in how we teach and who we teach\ldots and teaching play is hard work}" \cite{Su:2017}. Granting students opportunities to take what can sometimes be digestible, but dry material, such as interpolation, and allowing them to get their hands dirty by experiencing its utility in mathematical models at the interface of education and contemporary research, could have a profound impact on their future mathematical or scientific journeys.

%
%

\section*{Acknowledgments}
The author would like to thank Charles Peskin for the development of immersed boundary method and Boyce Griffith for IBAMR, to which many of the input files structures of \textit{IB2d} are based. He would also like to thank Austin Baird, Aaron Barrett, Christina Battista, Robert Booth, Karen Clark, Jana Gevertz, Christina Hamlet, Alexander Hoover, Shannon Jones, Andrea Lane, Laura Miller, Matthew Mizuhara, Arvind Santhanakrishnan, Michael Senter, Christopher Strickland, and Lindsay Waldrop for comments on the design of the \textit{IB2d} software and suggestions for examples. He also wants to acknowledge his Spring 2018 Numerical Analysis class at The College of New Jersey (Yaseen Ayuby, Shalini Basu, Gina Lee Celia, Rebecca Conn, Alexander Cretella, Robert Dunphy, Alyssa Farrell, Sarah Jennings, Edward Kennedy, Nicole Krysa, Aidan Lalley, Jason Miles, Jessica Patterson, Brittany Reedman, Angelina Sepita, Nicole Smallze, Briana Vieira, and Ursula Widocki) as the original motivators for this project. This project was funded by the NSF OAC-1828163 and the TCNJ Support of Scholarly Activity (SOSA) Grant, the Department of Mathematics and Statistics, and the School of Science at TCNJ.

%
%
\appendix
%
%

\section{Details regarding \textit{IB2d} and the Immersed Boundary Method (IB)}
\label{IB:Appendix}

Here we will touch upon the major points regarding the fluid-structure interaction software used for computations, \textit{IB2d}, as well as the numerical method it is built upon, the \textit{immersed boundary method} (IB).

%
%

\subsection{\textit{IB2d}}
\label{app:IB2d_info}

Biological fluid dynamics is a vast subject, in which nearly encompasses the entire natural world around us. From the way birds fly, fish swim, or the way you've taken a couple breaths in the past few seconds, fluid dynamics, or more precisely, fluid-structure interactions are ever present. Unfortunately, for such a significant practical area of mathematical modeling, it traditionally comes with a very steep learning curve, making it challenging to teach educational modules or give students meaningful first hand experience in course projects. Our open source software, \textit{IB2d}, was designed specifically for these purposes. It has two full implementations in high-level programming environments most familiar to most undergraduate students, MATLAB and Python. 

\textit{IB2d} was created to be used for both teaching and research purposes. It comes equipped with over $60$ built in examples that allow students to explore the world of fluid dynamics and fluid-structure interaction, from examples that illustrate fluid dynamics principles, such as flow around a cylinder for multiple Reynolds Numbers or the Rayleigh-Taylor Instability, to examples that purely illustrate interactions of a fluid with different immersed structure material properties to biological examples, such as jellyfish locomotion or embryonic heart development. Some of these examples are highlighted in \cite{Battista:2015,BattistaIB2d:2016,BattistaIB2d:2017}. Therefore \textit{IB2d} can be used for either course projects or homework assignments for a multitude of courses, ranging from mathematical modeling and mathematical biology courses to fluid mechanics to scientific computing. It has also been used for research purposes \cite{Pallasdies:2019,Miles:2019b}.

For these reasons, there have been tutorial videos created to help acquaint one with the software. All tutorial videos be found at \url{github.com/nickabattista/IB2d}:
\begin{itemize}
    %
    \item \textbf{Tutorial 1}: \url{https://youtu.be/PJyQA0vwbgU} \\
    \textit{An introduction to the immersed boundary method, fiber models, open source IB software, IB2d, and some FSI examples!}
    %
    \item \textbf{Tutorial 2}: \url{https://youtu.be/jSwCKq0v84s} \\
    \textit{A tour of what comes with the IB2d software, how to download it, what Example sub-folders contain and what input files are necessary to run a simulation}
    %
    \item \textbf{Tutorial 3}: \url{https://youtu.be/I3TLpyEBXfE} \\
    \textit{The basics of constructing immersed boundary geometries, printing the appropriate input file formats, and going through these for the oscillating rubberband example from Tutorial 2}
    %
    \item \textbf{Tutorial 4}: \url{https://youtu.be/4D4ruXbeCiQ} \\
    \textit{The basics of visualizing data using open source visualization software called VisIt (by Lawrence Livermore National Labs). Using the oscillating rubberband from Tutorial 2 as an example to visualize the Lagrangian Points and  Eulerian Data (colormaps for scalar data and vector fields for fluid velocity vectors)}
\end{itemize}

More explicit details about \textit{IB2d}'s functionality can be found in \cite{Battista:2015,BattistaIB2d:2016,BattistaIB2d:2017}.

%
%

\subsection{Governing Equations of IB}
\label{app:IB_eqns}

In this section we will introduce the equations of fluid motion and how they can be coupled with the motion and deformations of an immersed body. The conservation of momentum equations that govern an incompressible and viscous fluid are written as the following set of coupled partial differential equations,

\begin{equation} 
   \rho\Big[\frac{\partial\textbf{u}}{\partial t}({\textbf{x} },t) +\textbf{u}({\textbf {x}},t)\cdot\nabla \textbf{u}({\textbf{x}} ,t)\Big]=  -\nabla p({\textbf{x}},t) + \mu \Delta \textbf{u}({\textbf{x}},t) + \textbf{F}({\textbf{x}},t) \label{eq:NS1}
\end{equation}
\begin{equation}
      \nabla\cdot \textbf{u}({\bf x},t) = 0 \label{eq:NSDiv1}
\end{equation}
where $\textbf{u}({\textbf{x}},t) $ is the fluid velocity, $p({\textbf{x}},t) $ is the pressure, $\textbf{F}({\textbf{x}},t) $ is the force per unit area applied to the fluid by the immersed boundary, $\rho$ and $\mu$ are the fluid's density and dynamic viscosity, respectively. The independent variables are the time $t$ and the position ${\textbf{x}}$. The variables $\textbf{u}, p$, and $\textbf{F}$ are all written in an Eulerian frame on the fixed Cartesian mesh, $\textbf{x}$. We note that (\ref{eq:NS1}) is the conversation of momentum, while (\ref{eq:NSDiv1}) is the conversation of mass, for an incompressible fluid.

The equations that couple the motion of the fluid to deformations of the structure are written as integral equations. These \textit{interaction} equations handle all communication between the fluid (Eulerian) grid and immersed boundary (Lagrangian grid). They are given as the following integral equations with delta function kernels,
\begin{align}
   {\bf F}({\bf x},t) &= \int {\bf f}(s,t)  \delta\left({\bf x} - {\bf X}(s,t)\right) ds \label{eq:force1} \\
   {\bf U}({\bf X}(s,t))  &= \int \textbf{u}({\bf x},t)  \delta\left({\bf x} - {\bf X}(s,t)\right) d{\bf x} \label{eq:force2}
\end{align}
where ${\bf f}(s,t)$ is the force per unit length applied by the boundary to the fluid as a function of Lagrangian position, $s$, and time, $t$, $\delta({\bf x})$ is a three-dimensional delta function, and ${\bf X}(s,t)$ gives the Cartesian coordinates at time $t$ of the material point labeled by the Lagrangian parameter, $s$. The Lagrangian forcing term, ${\bf f}(s,t)$, gives the deformation forces along the boundary at the Lagrangian parameter, $s$. (\ref{eq:force1}) applies this force from the immersed boundary to the fluid through the external forcing term in (\ref{eq:NS1}). Equation (\ref{eq:force2}) moves the boundary at the local fluid velocity. This enforces the no-slip condition. Each integral transformation uses a three-dimensional Dirac delta function kernel, $\delta$, to convert Lagrangian variables to Eulerian variables and vice versa.

The way deformation forces are computed, e.g., the forcing term, $\textbf{f}(s,t)$, in the integrand of (\ref{eq:force1}), is specific to the application. To either hold the geometry nearly rigid or prescribe the motion of the immersed structure, all of the Lagrangian points along the immersed boundary are tethered to target points. They can do this through a penalty force formulation of $\textbf{f}(s,t)$. In this paper, in Sections \ref{Interpolation} and Section \ref{Heart_Pumping_Example}, we have used target points to prescribe the motion of the immersed structure. The penalty force was written in the following way,
\begin{equation}
{\bf f}(s,t) = k_{targ} \left(\textbf{Y}(s,t) - {\bf X}(s,t)\right),
\label{eq:force3}
\end{equation}
where $k_{targ}$ is a stiffness coefficient and $\textbf{Y}(s,t)$ is the prescribed position of the target boundary. Note that $\textbf{Y}(s,t)$ is a function of both the Lagrangian parameter, $s$, and time, $t$, and that in these models $k_{targ}$ was chosen to be large so that it would effectively drag the Lagrangian points into the preferred positions.

In Section \ref{sec:swimmer}, we construct a swimmer that is composed of springs and beams. Springs allow for stretching and compressing of the successive Lagrangian points, while beams allow for bending. Their corresponding deformation force equations can be written as the following,
\begin{align}
    \label{fiber_spring} \mathbf{F}_{spr} &= -k_{spr} \left( 1 - \frac{R_L}{\left|\left| \mathbf{X}_{S} - \mathbf{X}_M \right|\right| } \right) \cdot \left( \mathbf{X}_M - \mathbf{X}_S \right). \\
    \label{fiber_beam} \mathbf{F}_{beam} & =-k_{beam} \frac{\partial^4}{\partial s^4}\Big( \mathbf{X}(s,t) - \mathbf{X}_B(s,t) \Big),
\end{align}
where $k_{spr}$ and $k_{beam}$ are the spring stiffness and beam stiffness coefficients for springs and beams, respectively. For the linear spring forces, the terms $X_{M}$ and $X_{S}$ represent the positions in Cartesian coordinates of the master and slave Lagrangian nodes at time, $t$, and $R_L$ is the spring's corresponding resting length. For the bending force, $\mathbf{X}_B(s,t)$ represents the preferred curvature of the configuration at time, $t$. We note that in the swimmer model of Section \ref{sec:swimmer}, we interpolate between different curvature states given by different configurations of $\mathbf{X}^a_B(s,t)$ and $\mathbf{X}^b_B(s,t)$, rather than interpolate between positions in space for the swimmer.

Using delta functions as the kernel in (\ref{eq:force1})-(\ref{eq:force2}) is the heart of IB. To approximate these integrals, discretized (and regularized) delta functions are used. We use the ones given from \cite{Peskin:2002}, e.g., $\delta_h(\mathbf{x})$, 
\begin{equation}
\label{delta_h} \delta_h(\mathbf{x}) = \frac{1}{h^3} \phi\left(\frac{x}{h}\right) \phi\left(\frac{y}{h}\right) \phi\left(\frac{z}{h}\right) ,
\end{equation}
where $\phi(r)$ is defined as
\begin{equation}
\label{delta_phi} \phi(r) = \left\{ \begin{array}{l} \frac{1}{8}(3-2|r|+\sqrt{1+4|r|-4r^2} ), \ \ \ 0\leq |r| < 1 \\    
\frac{1}{8}(5-2|r|+\sqrt{-7+12|r|-4r^2}), 1\leq|r|<2 \\
0 \hspace{2.1in} 2\leq |r|.\\
\end{array}\right.
\end{equation}

%
%

\subsubsection{Numerical Algorithm}
\label{app:IB_Numerical_Algorithm}

As stated in the main text, we impose periodic and no slip boundary conditions on a rectangular domain. To solve \ref{eq:NS1}), (\ref{eq:NSDiv1}),(\ref{eq:force1}) and (\ref{eq:force2}) we need to update the velocity, pressure, position of the boundary, as well as the force acting on the boundary at time $n+1$ using data from time $n$. The IB does this in the following steps \cite{Peskin:2002,BattistaIB2d:2016}:

\textbf{Step 1:} Find the force density, ${\bf{F}}^{n}$ on the immersed boundary, from the current boundary configuration, ${\bf{X}}^{n}$.\\
\indent\textbf{Step 2:} Use (\ref{eq:force1}) to spread this boundary force from the Lagrangian boundary mesh to the Eulerian fluid lattice points.\\
\indent\textbf{Step 3:} Solve the Navier-Stokes equations, (\ref{eq:NS1}) and (\ref{eq:NSDiv1}), on the Eulerian grid. Upon doing so, we are updating ${\bf{u}}^{n+1}$ and $p^{n+1}$ from ${\bf{u}}^{n}$, $p^{n}$, and ${\bf{f}}^{n}$. \\
\indent\textbf{Step 4:} Update the material positions, ${\bf{X}}^{n+1}$, using the local fluid velocities, ${\bf{U}}^{n+1}$, computed from ${\bf{u}}^{n+1}$ and (\ref{eq:force2}).

\bibliographystyle{siamplain}
\bibliography{heart}
\end{document}


\maketitle

%
%

\section{Supplement 1:} The main supplementary file contains movies, codes, and slides pertaining to all the simulations detailed in this paper. It encompasses the following: $ $\\
\label{supp:1}

\begin{enumerate}
    \item Slides (both \textit{.pptx} and \textit{.pdf}) to streamline integration into a classroom. Note that movies are embedded in the \textit{.pptx} slides.
    \item Skeletal codes to run all the simulations (with the necessary source files, found in the \texttt{IBM$\_$Blackbox} folder, that can run upon downloading)
    \item The analysis scripts used to perform data analysis on the idealized anguilliform swimmer cases (with the necessary sources files, found in the\\ \texttt{IB2d$\_$Data$\_$Analysis$\_$Blackbox} folder, and simulation data, kept in \texttt{viz$\_$IB2d} folders, that can run upon downloading. Moreover all necessary paths are already set in the scripts.
    \item  Movies of all the simulations, to whose snapshots were shown in the manuscript. 
    \item VisIt Sessions \cite{HPV:VisIt} to visualize all the given data in the same manner that the snapshots in the manuscript did. For this purpose, upon opening VisIt one would go to \textit{File} $\rightarrow$\textit{Restore session with sources} option and then select the desired \texttt{VisIt$\_$Session.session} file and when prompted select the appropriate corresponding data. Note: the file paths in the VisIt session files cannot be automatically set in the software and must be adjusted upon use. Note that VisIt v2.12.3 was used for all the provided movies and snapshots.
\end{enumerate}

$ $ \\
We will now go into specifics about each sub-directory within this Supplemental Directory. $ $ \\

\begin{itemize}
    
    %
    %
    \item \texttt{interp$\_$Function$\_$Coeffs.m}: the script that is used to calculate the cubic interpolant coefficients.
    
    %
    %
    \item \textit{Circles/}: folder containing all simulation data from Section 2.
    \begin{enumerate}
        \item \texttt{Linear$\_$Interp}: contains codes to run the simulation and visualizing using VisIt, a movie of the corresponding simulation, and the data from having run said simulation. Note that this simulation was designed on a very coarse mesh so students could run it in approximately one minute on a personal machine.
        \item \texttt{Cubic$\_$Interp}: contains codes to run the simulation and visualizing using VisIt, a movie of the corresponding simulation, and the data from having run said simulation. Note that this simulation was designed on a very coarse mesh so students could run it in approximately one minute on a personal machine.
    \end{enumerate}

    %
    %
    \item \texttt{Pulsing$\_$Heart/}: folder containing all simulation data from Section 3, e.g., contains codes to run the simulation and visualizing using VisIt, a movie of the corresponding simulation, and the data from having run said simulation. Note that this simulation was designed on a very coarse mesh so students could run it in a couple of minutes on a personal machine.

    %
    %
    \item \texttt{Swimmer}: folder containing sub-folders that correspond to each case in Section 4. Note that these simulations were designed on a fine mesh and each takes on the order of $\sim2$ hours to run on a personal machine ($\sim$4-16GB RAM, $\sim$2-3GHz processor); the fine mesh is required for locomotion.
    \begin{enumerate}
        \item \texttt{Individual$\_$Swimmer}: the skeletal simulation code accompanied by Lagrangian data and select Eulerian data (in viz$\_$IB2d) that can analyzed using the provided analysis script, a movie of the swimmer, and its corresponding VisIt session file used to create the movie.
        \item \texttt{Case1}: for simulations when varying $(p_1,p_2)$ `symmetrically'; all skeletal simulation codes accompanied by Lagrangian data (in viz$\_$IB2d) that can analyzed using the included analysis script, movie of the swimmers, and corresponding VisIt session file used to create the movie.
        \item \texttt{Case2}: for simulations when varying $(p_1,p_2)$ `asymmetrically'; all skeletal simulation codes accompanied by Lagrangian data (in viz$\_$IB2d) that can analyzed using the included analysis script, movie of the swimmers, and corresponding VisIt session file used to create the movie.
        \item \texttt{Case3}: for simulations when varying the upstroke percentage of the total stroke period; all skeletal simulation codes accompanied by Lagrangian data (in viz$\_$IB2d) that can analyzed using the included analysis script, movie of the swimmers, and corresponding VisIt session file used to create the movie.
        \item \texttt{Viscosity$\_$Race}: a movie comparing the swimmers each immersed in a different fluid of different viscosity. Raw Lagrangian data was not included, but simulations can be run by changing the viscosity appropriately in the \textit{input2d} file.
        \item \texttt{IB2d$\_$Data$\_$Analysis$\_$Blackbox}: scripts used to perform the data analysis that do not need to be modified. Note these are the same Data Analysis scripts that \textit{IB2d} offers (as of August 23, 2018).
    \end{enumerate}
    
\end{itemize}

%
%

\section{Supplement 2: Spline Interpolant Coefficients for Swimmers}
\label{app:swimmer:coeffs}

In this supplemental section we list the spline interpolation coefficients when varying $(p_1,p_2)$ for our swimmer in Section 4 and solving the linear system from Section 2, (2.14). The script used to solve this linear system is given in Supplement 1 (\texttt{Supplemental/interp$\_$Function$\_$Coeffs.m}) as well as in the \textit{IB2d} directory: Examples/Examples$\_$Education/Interpolation/ from \url{https://github.com/nickabattista/ib2d}. We will list coefficients for both the symmetric and asymmetric cases. 

%
%
\subsection{Symmetric ($p_1,p_2$) coefficients} $ $\\

\begin{enumerate}
    \item $(p_1,p_2)=(0.1,0.9)$
    
    \begin{equation}
    \label{coeffs1} \begin{array}{lll} 
    a_0 = 0                & b_0 = 0.014   & c_0 = -10.111\\ 
    a_1 = 0                & b_1 = -0.417  & c_1 = 33.333 \\
    a_2 = 0                & b_2 = 4.167   & c_2 = -33.333 \\
    a_3 = 11.1111\ \ \ \ \ \ & b_3 = -2.778\ \ \ \ \ \ & c_3 = 11.111 \\\end{array}
    \end{equation}

    \item $(p_1,p_2)=(0.2,0.8)$
    
    \begin{equation}
    \label{coeffs2} \begin{array}{lll} 
    a_0 = 0                & b_0 = 0.083   & c_0 = -5.250\\ 
    a_1 = 0                & b_1 = -1.250  & c_1 = 18.750 \\
    a_2 = 0                & b_2 = 6.250   & c_2 = -18.750 \\
    a_3 = 6.250\ \ \ \ \ \ & b_3 = -4.167\ \ \ \ \ \ & c_3 = 6.250 \\\end{array}
    \end{equation}

    \item $(p_1,p_2)=(0.3,0.7)$
    
    \begin{equation}
    \label{coeffs3} \begin{array}{lll} 
    a_0 = 0                & b_0 = 0.321   & c_0 = -3.762\\ 
    a_1 = 0                & b_1 = -3.214  & c_1 = 14.286 \\
    a_2 = 0                & b_2 = 10.714   & c_2 = -14.286 \\
    a_3 = 4.762\ \ \ \ \ \ & b_3 = -7.143\ \ \ \ \ \ & c_3 = 4.762 \\\end{array}
    \end{equation}

    \item $(p_1,p_2)=(0.4,0.6)$

    \begin{equation}
    \label{coeffs4} \begin{array}{lll} 
    a_0 = 0                & b_0 = 1.333   & c_0 = -3.167\\ 
    a_1 = 0                & b_1 = -10.000  & c_1 = 12.500 \\
    a_2 = 0                & b_2 = 25.000   & c_2 = -12.500 \\
    a_3 = 4.167\ \ \ \ \ \ & b_3 = 16.667\ \ \ \ \ \ & c_3 = 4.167 \\\end{array}
    \end{equation}

\end{enumerate}

%
%
\subsection{Asymmetric ($p_1,p_2$) coefficients} $ $\\

\begin{enumerate}

    \item $(p_1,p_2)=(0.1,0.9)$
    
    \begin{equation}
    \label{coeffs5} \begin{array}{lll} 
    a_0 = 0                & b_0 = 0.014   & c_0 = -10.111\\ 
    a_1 = 0                & b_1 = -0.417  & c_1 = 33.333 \\
    a_2 = 0                & b_2 = 4.167   & c_2 = -33.333 \\
    a_3 = 11.111\ \ \ \ \ \ & b_3 = -2.778\ \ \ \ \ \ & c_3 = 11.111 \\\end{array}
    \end{equation}

    \item $(p_1,p_2)=(0.1,0.7)$
    
    \begin{equation}
    \label{coeffs6} \begin{array}{lll} 
    a_0 = 0                & b_0 = 0.019   & c_0 = -2.704\\ 
    a_1 = 0                & b_1 = -0.556  & c_1 = 11.111 \\
    a_2 = 0                & b_2 = 5.556   & c_2 = -11.111 \\
    a_3 = 14.286\ \ \ \ \ \ & b_3 = -4.233\ \ \ \ \ \ & c_3 = 3.704 \\\end{array}
    \end{equation}

    \item $(p_1,p_2)=(0.1,0.5)$
    
    \begin{equation}
    \label{coeffs7} \begin{array}{lll} 
    a_0 = 0                & b_0 = 0.028   & c_0 = -1.222\\ 
    a_1 = 0                & b_1 = -0.833  & c_1 = 6.667 \\
    a_2 = 0                & b_2 = 8.333   & c_2 = -6.667 \\
    a_3 = 20.0\ \ \ \ \ \ & b_3 = -7.778\ \ \ \ \ \ & c_3 = 2.222 \\\end{array}
    \end{equation}

    \item $(p_1,p_2)=(0.1,0.3)$

    \begin{equation}
    \label{coeffs8} \begin{array}{lll} 
    a_0 = 0                & b_0 = 0.056   & c_0 = -0.587\\ 
    a_1 = 0                & b_1 = -1.667  & c_1 = 4.762 \\
    a_2 = 0                & b_2 = 16.667   & c_2 = -4.762 \\
    a_3 = 33.333\ \ \ \ \ \ & b_3 = -22.222\ \ \ \ \ \ & c_3 = 1.587 \\\end{array}
    \end{equation}

\end{enumerate}

%
%

\bibliographystyle{siamplain}
\bibliography{heart}